\documentclass[12pt]{article}
\usepackage{amsmath}
\usepackage{amssymb}
\usepackage[]{graphicx}
\usepackage{epsfig, float}
\usepackage[cp1251]{inputenc}
\usepackage[english]{babel}
\usepackage{color}

\newtheorem{proposition}{Proposition}

\newtheorem{corollary}{Corollary}
\newcommand{\be}{\begin{equation}}
\newcommand{\ee}{\end{equation}}
\newcommand{\bea}{\begin{eqnarray}}
\newcommand{\eea}{\end{eqnarray}}

\newcommand{\ep}{\varepsilon}

\newcommand{\cop}{{\cal COP}^p}
\newcommand{\cp}{{\cal CP}^p}

\newcommand{\cb}{{\cal B}}
\newcommand{\ca}{{\cal A}}

\newcommand{\bftau}{{{\boldsymbol {\tau}}}}

\newtheorem{theorem}{Theorem}
\newtheorem{remark}{Remark}

\newtheorem{lemma}{Lemma}
\newtheorem{definition}{Definition}

 \newcounter{myc}

\begin{document}

 \date{}

\title{Explicit criterion of uniform  LP  duality for linear problems of copositive optimization
\thanks{This work was partially supported state research program "Convergence" (Republic Belarus), Task 1.3.01, by
 Portuguese funds through CIDMA - Center for Research and Development in Mathematics and Applications, and FCT -
 Portuguese Foundation for Science and Technology, within the project  UIDB/04106/2020.}}

\author{Kostyukova O.I.\thanks{Institute of Mathematics, National Academy of Sciences of Belarus, Surganov str. 11, 220072, Minsk,
 Belarus  ({\tt kostyukova@im.bas-net.by}).} \and Tchemisova T.V.\thanks{Mathematical Department, University of Aveiro, Campus
Universitario Santiago, 3810-193, Aveiro, Portugal ({\tt tatiana@mat.ua.pt}).}\and Dudina O.S.\thanks{Department of Applied Mathematics and Computer Science, Belarusian State University,  Nezavisimosti Ave., 4, 220030, Minsk, Belarus
({\tt dudina@bsu.by}).}}
\maketitle

\begin{abstract}
 An uniform LP duality is an useful property of conic matrix systems.   A consistent linear conic optimization  problem
yields  uniform LP duality if for any linear cost function, for which  the primal problem has finite optimal value, the corresponding Lagrange dual problem is attainable
and the duality gap vanishes.

 In this paper,  we establish  new  necessary and sufficient conditions guaranteing  the  uniform LP duality
for linear problems of Copositive Programming and   formulate these conditions in different equivalent
forms. The main results are obtained   using  an
 approach  developed  in previous papers of the authors   and based  on a concept of immobile indices
that permits  alternative  representations of the {set} of feasible solutions.

\end{abstract}

\section{Introduction}

Conic optimization is a subfield of convex optimization devoted to problems of minimizing a
 convex function over an intersection of an affine subspace and a convex cone.
 Conic  problems form a broad and important class of optimization problems since, according to \cite{Nesterov},
any convex optimization problem can be represented as a conic  one.
 This class  includes some of the most well-known types of convex  problems,
such as linear and semidefinite programming  problems ((LP) and (SDP), respectively).  Many
 problems of semi-infinite programming (SIP) consisting in optimization of a cost function w.r.t.
 an infinite number of functional constraints, can be also considered as conic optimization problems.

Copositive Programming (CoP)  can be thought of as   a special  case of SIP and   a generalization of SDP.
	  In CoP,   a linear function is optimized over a cone of  matrices
	that are positive semidefinite in the non-negative ortant $\mathbb{R}^p_+$  ({\it copositive matrices}).
	 Formally, problems of CoP are very similar to that of SDP, but   CoP
	deals with  more complex and less studied  problems than  SDP.  Being a fairly new field of
	research, CoP  has already gained popularity, 	as it   has  been   proven
	 to be  very useful in modeling particularly  complex  problems of
	 convex    optimization, graph theory, algebra, and different applications
	(see, for example, \cite{Bomze_2012,Dur2010}, and the references there).

 Optimality conditions and   the  duality relations are among the most
emerging optimization topics, and the importance of studying them has long been recognized (see e.g. \cite{borwein,Duffin,KTD,Li,Ramana-W}, and the references therein).
Duality plays a central role in testing optimality, identifying infeasibility, establishing lower-bounds of optimal objective value, design and analysis of iterative
algorithms.

Traditionally, for a given (primal) convex problem, based on its initial data,
the Lagrangian dual problem   is constructed. The difference between optimal values of primal and dual cost functions is called the {\it duality gap}.

 The primal and the (Lagrange)  dual problems are closely related. The strength of this relationship depends on the initial problem data, which specify the constraints and the
 cost function.
Roughly speaking, a pair of dual problems is said to satisfy
(i) a {\it weak duality} if the duality gap is non-negative, (ii) a {\it strong duality} if, for a given cost function, the duality gap is zero, and (iii) an {\it uniform
duality} if  the  duality gap is zero for any cost function.

 It is well-known that  the strong duality is guaranteed  unconditionally
only for  the  LP problems, while for most important classes of convex and
conic problems this property is satisfied only under certain rather strong assumptions.
 Many  works are devoted to study of these assumptions (see e.g. \cite{borwein1992,Fang,Jeyak}, and the references therein).

  For a linear SIP problem,  one of the  known   criteria
	of the uniform duality is  the   closedness of  some   cone   built  on the basis
	of  the  problem's  data (see \cite{Duffin}).
	The criterion can be applied to any linear SIP problems, but it has an implicit form. That is why   researchers often
	try to find more explicit conditions for  the  uniform duality, taking into account  the  specifics
	of the  problems under consideration.

The concept of a {\it uniform LP duality}  was   considered  in the work  of Duffin et.al (\cite{Duffin})  for
linear SIP problems, and in \cite{Jeroslow}, it has been  used  for a wider class of convex SIP problems in the form of an    {\it uniform convex  duality}.  In   {\cite{Li,
Pataki2017, Zhang}},   the uniform duality
property was studied for certain classes of linear conic, convex
SIP and SDP problems and the conditions guaranteing that this property is   satisfied, were deduced.

Not much literature is available for optimality  and duality conditions for CoP problems.
Moreover, the strong/uniform duality in CoP is not easy to establish due to  intrinsic complexity
of copositive problems connected with the fact that the  cone  of copositive matrices and the corresponding dual cone
of {\it completely positive} matrices   {do} not possess some  "good"  properties:
 they are neither  self-dual, no facially exposed, no symmetric. At the same time, a more in-depth study of the duality  issues  and the description of explicit criteria for
 the fulfillment of the properties of strong/uniform duality is an extremely important challenge not only for the theory of the CoP, but also for the development of algorithms
 and numerical applications.

 The aim of this paper is to establish  new  necessary and sufficient conditions guaranteing  the  uniform LP duality
for linear CoP problems, and to formulate these conditions in different equivalent
forms thus broadening their scope. The main results are obtained on the base of an
  approach  developed  in previous papers of the authors.
This  approach is based on a concept of immobile indices and  first  was  {described} for
SIP problems (see, for example,    \cite{JOTA10}),  and then applied to various classes of convex conic problems in
 \cite{KTD,KT-IM}, and others. The concept of immobile indices and  {the} properties of
 the sets generated by these indices, make  it possible to constructively represent for the CoP problem some important
 subcones  used in conic optimization  (the  faces of the copositive cone, in particular, the minimal  face containing the feasible set)
  to  obtain new CQ-free  strong dual formulations and  optimality  conditions.
	  In this paper,     the   set of immobile indices  is used  to obtain  new  criteria  of the  uniform duality for linear CoP problems.

The remainder of this paper is organized as follows. The problem's statement and relevant research is overviewed in section 2.
 The main results of the paper, new necessary and sufficient conditions of uniform duality for linear copositive problems, are formulated and proved in section 3.
Several equivalent formulations of the  uniform duality conditions from section 3, are deduced in section 4.
 Section 5 contains examples that confirm that the conditions obtained in the paper are essential.   Some comparison with known results is given.
 In section 6, we analyze  the  uniform duality conditions for SIP problems applied to CoP. We show that
  results  obtained in this paper   allow one to reformulate these conditions in a more   explicit  form.
The final section 7 contains some conclusions, and several technical proofs are situated in Appendix.

\section{Problem statement and preliminary results}\label{S-2}

Given a finite-dimensional vector space $\mathfrak{X}$,
let's, first, recall  some generally accepted definitions.

\vspace{2mm}

A  set $C\subset \mathfrak{X}$ is {\it convex} if for any $x, y \in C$ and any $\alpha \in [0,1]$, it holds $\alpha x+(1-\alpha ) y \in C.$
 A set $K \subset \mathfrak{X}$ is a {\it cone} if for any $x \in K$ and any $\alpha > 0$, it holds $\alpha x\in K.$
Given a cone $K \subset \mathfrak{X}$, its {\it dual cone} $K^*$
is given by \begin{equation*}K^*:=\{x\in \mathfrak{X}: \langle x,y\rangle \geq 0  \ \forall y\in K\}.\end{equation*}

Given a set $\mathcal B \mathfrak{X},$  denote by
${\rm conv } \mathcal B$   its {\it convex hull}, i.e., the minimal (by inclusion)
convex set, containing this set, by ${\rm span}(\mathcal B)$    its {\it span}, i.e.,
 the smallest linear subspace  containing  $\mathcal B, $  and  by ${\rm cone } {\cal B}$ its
  {\it conic hull}, i.e the set of all conic combinations of the points   of  ${\cal B}$.
 In what follows, we will denote by ${\rm cl}\, {\cal B}$  the {\it closure} of the set ${\cal B}$, by  $ {\rm int} \, {\cal B}$ its {\it   interior}, and by  $ {\rm relint}\, {\cal B}$ its {\it relative interior}.

 Given an integer $p>1$, consider  the vector space $\mathbb R^p$ with the standard orthogonal basis
 $\{{\bf e}_k,  \ k=1,2,\dots, p\}$. Denote by  $\mathbb R^{p}_+$  the set of all $ p$ - vectors  with non-negative components,
   by ${\mathcal S}^{p}$    the space of real symmetric $p\times p$ matrices,   and by ${\cal S}^p_+$ the cone of symmetric  positive semidefinite $p\times p$ matrices. The
   space ${\mathcal S}^{p}$ is considered here as a vector space with the {\it trace inner product} $A\bullet B:={\rm trace}\, (AB).$

In this paper, we  deal  with  special  classes of    cones, the elements   of  which are   matrices,
 in particular,  with the cones of   {\it copositive} and {\it completely positive} matrices.

 Let $\mathcal{COP}^{p}$  denote the cone of symmetric copositive $p\times p$ matrices:
\begin{equation*}\mathcal{COP}^p:=\{D\in {\mathcal S}^{p}:{{\bf t}}^{\top}D{{\bf t}}\geq 0\ \forall {\bf t} \in \mathbb
R^p_+\}.\end{equation*}
 Consider a compact subset of
 $\mathbb R^{p}_+$ in the form of the simplex \begin{equation}\label{SetT}T:=\{{\bf t} \in \mathbb R^p_+:\mathbf{e}^{\top}{\bf t}=1\} \end{equation}
with $\mathbf{e}=(1,1,...,1)^{\top}\in \mathbb R^p $.
It is evident that the cone $\cop$ can be equivalently described as follows:
\begin{equation}\cop=\{D\in \mathcal S^{p}: \ {\bf t}^{\top}D{\bf t}\geq 0 \ \forall{\bf t}\in T\}.\label{copT}\end{equation}
 The dual cone to $\cop$  is the cone of   completely positive matrices defined as
\begin{equation*}{(\cop)^*=\mathcal{CP}^p:}= {\rm conv}\{ {\bf t}{\bf t} ^{\top}: \  {\bf t} \in {\mathbb R}^p_+\}.\end{equation*}

 The  cones of  copositive and completely positive  matrices  are  known to be  {\it proper}   cones, which  means that   they are  closed, convex, pointed, and  full-dimensional.

Consider a linear copositive programming  problem in the form
$${ \mbox{ \bf P:}}\;\;\qquad  \qquad  \qquad \qquad  \displaystyle\min_{{\bf x}  \in \mathbb R^n } \ {\bf c}^{\top}{\bf x} \;\;\;
\mbox{s.t. } {\cal A}({\bf x})
\in \mathcal{COP}^p,\qquad \qquad \qquad \qquad  $$
 where $ {\bf x}=(x_1,...,x_n)^{\top}$  is the vector of  decision variables,  the constraint matrix function $\mathcal A({\bf x})$ is defined as
$ \mathcal A({\bf x}):=\displaystyle\sum_{m=1}^{n}A_mx_m+A_0;$
vector $ {\bf c}\in \mathbb R^n $  and matrices $A_m\in \mathcal S^{p},$ $i=0,1,\dots,n$  are given.  Denote by $X$  the set of feasible solutions of this problem:
$$X=\{{\bf x} \in \mathbb R^n:{\cal A}({\bf x})\in \mathcal{COP}^p\}.$$

 For  the  problem ({\bf P}),  the {\it Lagrange dual problem} takes the form:
$$ \mbox{\bf D}:\ \ \qquad  \max  -U\bullet A_0, \mbox{ s.t. } U\bullet A_m=c_m \;  \forall  m=1,2,...,n;\ U\in \cp.\qquad $$

In what follows, for an optimization problem ({\bf Q}), $Val({\bf Q})$ denotes the optimal value of the objective function in the problem
 ({\bf Q}) (shortly, the optimal value of the problem ({\bf Q})).

 \vspace{2mm}

It is a known fact  (see, for example,  \cite{KT_strong}  and section \ref{s:5} below)  that  for CoP   problems,  the
  optimal values $Val ({\bf P} )$ and $Val( {\bf D} )$ of the primal problem $({\bf P}) $
and the corresponding Lagrange dual problem $(\bf D)$ are not necessarily equal, even if they exist and are finite.
 A situation where, assuming   $Val ({\bf P})>-\infty$,  the problem   $(\bf D)$ has
 optimal solution and the so-called {\it duality gap}, the difference $Val ({\bf P} )-Val( {\bf D} )$
   equals to zero, is called  a {\it strong duality}. The paper \cite{KT_strong}  is  devoted to strong dual  formulations for copositive problems that   differ  from  the
   Lagrange dual problem ({\bf D}).
   
    In this paper,  for linear CoP problems,  we consider  a slightly different
 duality property of their feasible sets that guarantees vanishing of the duality gap for all   cost vectors ${\bf c}$.
 Since this property is related to   {the constraint system  $\ca({\bf x})\in \cop$ of  the problem $({\bf P}) $,} we will refer to
  this property as to a property of  {this}  system.

\begin{definition}\label{def1}  A consistent system $\ca({\bf x})\in \cop$ yields  the property of {\it uniform LP duality}  if for all ${\bf c}\in \mathbb R^n$, such that the
optimal  value of the problem  $(\bf P)$ is finite ($Val({\bf P})>-\infty $),  the corresponding Lagrange dual  problem ({\bf D}) has an optimal solution, and  it holds
$Val ({\bf P})=Val({\bf D}).$
\end{definition}
  It is known that under  the  Slater condition (the property that   for some ${\bf x}\in \mathbb R^n$, it holds $\ca({\bf x})\in {\rm int} (\cop)$ ), the   system $\ca({\bf x})\in \cop$ yields the  uniform LP duality property.

\vspace{2mm}

 Given  the set  $X$  of feasible solutions of the problem (\textbf{P}), denote  by $T_{im}$  the  set of {\it normalized immobile indices} of constraints  in this problem:
$$T_{im}:=\{{\bf t} \in T: {\bf t}^\top\ca({\bf x}){\bf t}=0\ \forall {\bf x} \in X\}.$$

   Some  properties of the set $T_{im}$  in copositive problems
	 were established in our previous works (see e.g. \cite{KTD, KT-new,KT-Math}).
In particular, it is known  that the set $T_{im}$  is either empty or an union of a finite number of convex bounded polyhedra.
 Also, it was shown in \cite{KTD} that the emptiness of the set $T_{im}$ is equivalent to the fulfillment of the Slater condition.

Suppose that  the set  $T_{im}$   is not empty  and denote by $${\cal T}:=\{{\bftau}(j), j \in  J\}, \ |J|<\infty,$$ the set of vertices of   its convex hull  ${\rm
conv}T_{im}.$
  It was shown in \cite{KTD} that
\be  X\subset Z:=\{{\bf x} \in \mathbb R^n: \ca({\bf x}){\bftau}(j)\geq 0,\ j \in J\}.\label{06-12-1}\ee

Denote $P:=\{1,...,p\} $  and  introduce the following sets:
\be \overline M(j):=\{k \in P:  {\bf e}^\top_k\ca({\bf x}){\bftau}(j)=0 \ \forall {\bf x} \in Z\}, \ j \in J,\label{Lj}\ee
\be M(j)=\{k \in P:  {\bf e}^\top_k\ca({\bf x}){\bftau}(j)=0 \ \forall {\bf x} \in X\},\  j \in J,\label{Mj}\ee
 \be N_*(j):=\{k \in P: \ \exists \,  {\bf x}(k,j)\in X \mbox{ such that }  {\bf e}^\top_k\ca({\bf x}(k,j)){\bftau}(j)=0\}, \ j \in J.\label{W5}\ee

The following lemma is proved in \cite{KT-IM}.
\begin{lemma}\label{KOST-L4}  Given a  consistent system $\ca({\bf x})\in \cop$  with
 the corresponding set of  normalized immobile indices of constraints  $T_{im}$ and   the sets $\overline M(j),$ $M(j)$, $j \in J,$
 defined in (\ref{Lj}) and (\ref{Mj}) respectively, it holds
$ \overline M(j)=M(j) \ \forall j \in J.$
\end{lemma}

It was shown in \cite{KT-new} that  the  problem $({\bf P})$ is equivalent to the following problem:
 \bea \mbox{{\bf P}}_{*}:\ \ \qquad &  \qquad  \qquad \qquad  \qquad  \min {\bf c}^\top {\bf x} \qquad  \qquad  \qquad  \qquad \qquad  \qquad  \nonumber\\
 &\mbox{ s.t. } {\bf t}^\top \ca({\bf x}){\bf t}\geq 0 \ \forall {\bf t} \in \Omega,\
 {\bf e}^\top_k\ca({\bf x}){\bftau}(j)
 \geq 0 \ \forall k \in N_{*}(j)\  \forall j \in J, \ \ \label{W6}\eea
and there exists a {\it  minimally active} feasible solution  ${\bf x}^*\in X$ such that
\bea & {\bf t}^\top \ca({\bf x}^*){\bf t}>0  \; \forall {\bf t} \in T\setminus T_{im},\label{W7}\\
 &{\bf e}^\top_k\ca( {\bf x}^*){\bftau}(j)= 0,\, \forall k \in M(j);\;  {\bf e}^\top_k\ca({\bf x}^*){\bftau}(j)> 0 \, \forall k \in P\setminus  M(j), \forall j \in
 J.\nonumber\eea
Here and in what follows, we will use the  set
\be  \Omega:=\{{\bf t} \in T: \rho({\bf t}, {\rm conv} T_{im})\geq \sigma\}\subset T\setminus T_{im}, \ \label{Omega}\ee
where  $ \sigma:\!=\!\min\{\tau_k(j), k \in P_+({\bftau}(j)), j \in J\}>0, $
$P_+({\bf t}):=\{k \in P:t_k>0\}$ for ${\bf t}=(t_k,k \in P)^\top $ $\in \mathbb R^p_+$,
 $\ \rho({\bf t}, {\cal B} )=\min\limits_{{\bftau}\in  {\cal B}}\sum\limits_{k\in P}|t_k-\tau_k|$ for some set ${\cal B}\subset \mathbb R^p$.
\begin{remark} In  \cite{KT-new}, it was formulated a similar problem equivalent to the  problem $({\bf P})$.  This problem has  additional constraints
${\bf e}^\top_k\ca({\bf x}){\bftau}(j) \geq 0$ for all $k\in P\setminus (N_{*}(j)),$ $ j \in J.$
It worth mentioning that  these constraints can be omitted  since    they are non-active for all ${\bf x}\in X.$
\end{remark}

For $ {\bf t} \in T$ and $k\in P$, denote
\be \bar {\bf b}(k,{\bf t})=\left(\begin{array}{c}
{\bf e}^\top_kA_m{\bf t}\cr
m=0,1,...,n\end{array}\right)\in \mathbb R^{n+1},{\bf a}({\bf t})=\left(\begin{array}{c}
{\bf t}^\top A_m{\bf t}\cr
m=0,1,...,n\end{array}\right)=\sum\limits_{k\in P_+({\bf t})}t_k\bar{\bf b}(k,{\bf t})\in \mathbb R^{ n+1},\label{a-t}\ee
\be {\bf b}(k,j)=\left(\begin{array}{c}
{\bf e}^\top_kA_m{\bftau}(j)\cr
m=0,1,...,n\end{array}\right)=\bar {{\bf b}}(k,{\bftau}(j))\in \mathbb R^{n+1}.\label{b-t}\ee
\begin{proposition}\label{P0}  In the terminology of Lemma \ref{KOST-L4}, for any $j_0\in J$ and $k_0\in M(j_0),$ there exist numbers $\alpha_{kj}=
\alpha_{kj}(k_0,j_0),$ $ k\in M(j),\ j \in J,$ such that
\be -{\bf b}(k_0,j_0)=\sum\limits_{j\in J}\sum\limits_{k\in M(j)}\alpha_{kj}{\bf b}(k,j),\ \alpha_{kj}\geq 0\ \forall k\in M(j),\ j \in J,\label{07-12-1}\ee
 where vectors ${\bf b}(k,j)$ are defined in (\ref{b-t}).
\end{proposition}
{\bf Proof.}  Using  the  notation introduced above, the set $Z$   defined in (\ref{06-12-1})  can be represented as
\be Z=\{{\bf x} \in \mathbb R^n: (1, {\bf x}^\top){\bf b}(k,j)\geq 0\ \forall  k\in P,\,  \forall  j \in J\}.\label{07-12-4}\ee
Consider the following LP problem:
$$ {{\bf LP}_*}:\qquad\qquad \qquad \max (1, {\bf z}^\top){\bf b}(k_0,j_0)\quad  {\rm s.t. }\  {\bf z} \in  Z.\qquad \qquad $$ 
Due to Lemma \ref{KOST-L4}, we have  $k_0\in M(j_0)=\overline M(j_0)$ and it follows from the definition of the set $\overline M(j_0)$ that $Val({\bf LP}_*)=0.$
Hence a vector $ {\bf x}^*$ satisfying (\ref{W7}) is an optimal solution of the problem (${\bf LP}_*$). Taking into account representation
(\ref{07-12-4}) of the set $Z$ and  relations (\ref{W7}) we see that  relations  (\ref{07-12-1})  are necessary and sufficient
  optimality conditions  for ${\bf x}^*$ in  the problem (${\bf LP}_*$).  $\ \Box$

  \vspace{3mm}

  Let us partition the set $J$  into subsets $J(s),$ $s \in S$, as it was done in \cite{KT-Math}.
Then (see \cite{KT-Math})
\be T_{im}=\bigcup\limits_{s\in S} T_{im}(s), \  T_{im}(s):={\rm conv} \{{\bftau}(j), j \in J(s)\}, s \in S,\label{W32}\ee
\be P_*(s):=\bigcup\limits_{j \in J(s)}P_+({\bftau}(j))\subset M(j) \  \forall j \in J(s),   \ \forall  s \in S.\label{07-12-5}\ee

\begin{proposition}\label{Pr-07-12-1}  Consider  a consistent system $\ca({\bf x})\in \cop$
 with the corresponding sets  $T_{im}, $ ${\cal T}, $ $M(j),\ j \in J$,  and the vectors ${\bf a}({\bf t}), {\bf b}(k,j), k\in M(j),\ j \in J$ defined above. The following
 inclusions hold true:
\be {\bf a}({\bf t})\in {\rm cone}\{{\bf b}(k,j),\  k\in M(j),\ j \in J\} \ \ \forall  {\bf t}\in T_{im}.\label{07-12-6}\ee

\end{proposition}
{\bf Proof.} Consider any ${\bf t}\in T_{im}.$ It follows from (\ref{W32}) that ${\bf t}\in T_{im}(s)$ with some $s\in S.$ Hence
\be {\bf t}=\sum\limits_{j\in J(s)}\alpha_j\bftau(j),\ \alpha_j\geq 0 \ \forall  j \in J(s);\ \sum\limits_{j\in J(s)}\alpha_j=1.\label{W20}\ee
 Then we obtain
 \be {\bf a}({\bf t})=\sum\limits_{k\in P_+({\bf t})}t_k\bar{\bf b}(k,{\bf t})=\sum\limits_{k\in P_+({\bf t})}t_k\sum\limits_{j \in J(s)}\alpha_j{\bf b}(k,j)
=\sum\limits_{j \in J(s)}\sum\limits_{k\in M(j)}t_k\alpha_j{\bf b}(k,j).\label{07-12-7}\ee
Here we took into account that $P_+({\bf t})\subset \bigcup\limits_{j\in J(s)}P_+(\bftau^*(j))=P_*(s)$ and (\ref{07-12-5}).

Since $t_k\alpha_j\geq 0$ for all $k\in M(j)$ and $j\in J(s)\subset J,$ we conclude  from  (\ref{07-12-7}) that inclusions (\ref{07-12-6}) take place.
 $\  \Box$

\section{ Necessary and sufficient uniform  LP  duality conditions}\label{s:3}
In this section, we will prove two statements   containing new necessary and sufficient uniform  LP  duality conditions for  linear   CoP   systems.

\begin{proposition}\label{P1} A  consistent  linear CoP
 system $\ca({\bf x})\in \cop$ yields the  uniform LP duality iff the following relations hold:
\be  {\bf b}(k,j)\in {\rm cone}\{{\bf a}({\bf t}), {\bf t}\in T\}   \ \forall    k\in N_{*}(j), \    \ \forall j \in J.\label{cond-2}\ee
\end{proposition}
{\bf Proof.}   Notice that if $T_{im}=\emptyset$, then $J=\emptyset$, and we consider that conditions   (\ref{cond-2}) are fulfilled.

$\Rightarrow)$ Suppose that  the   consistent system $\ca({\bf x})\in \cop$ yields  the  uniform LP duality.  Then   for any ${\bf c}\in \mathbb R^n$   for which  $Val({\bf
P})>-\infty$, there exists a matrix $U=U({\bf c})$ in the form
\be U=\sum\limits_{i\in I}\alpha_i{\bf t}(i)({\bf t}(i))^\top, \ \alpha_i>0,\ {\bf t}(i)\in T,\ i\in I,\ |I| <\infty,\label{Uform}\ee
such that
\be A_m\bullet U=c_m,\, m=1,...,n; \  A_0\bullet U=-Val({\bf P}).\label{W1}\ee

For   fixed $j\in J$ and $k\in N_{*}(j)$, consider   the  problem ({\bf P}) with
$${\bf c}^\top=(c_m= {\bf e}^\top_k A_{m}{\bftau}(j), m=1,...,n).$$ It follows from  (\ref{06-12-1}) that
$$ {\bf c}^\top {\bf x}=\sum\limits_{m=1}^n {\bf e}^\top_k A_m x_{m}{\bftau}(j)\geq - {\bf e}^\top_k A_0 {\bftau}(j)\ \ \forall  {\bf x}\in X, $$ 
and it follows from (\ref{W5}) that there exists ${\bf x}(k,j)\in X$ such ${\bf c}^\top  {\bf x}(k,j)=- {\bf e}^\top_k A_0{\bftau}(j)$. Thus we can conclude that
  $Val({\bf P})=-{\bf e}^\top_k A_0 \bftau(j) >-\infty$. Taking into account that  the  system $\ca({\bf x})\in \cop$ yields  the  uniform  LP  duality, we conclude that
  there exists a matrix $U$ in the form (\ref{Uform}) such that
\bea & A_m\bullet U=\sum\limits_{i\in I}\alpha_i({\bf t}(i))^\top A_m{\bf t}(i)={\bf e}^\top_k A_{m}{\bftau}(j),\, m=1,...,n;\nonumber\\
& A_0\bullet U=\sum\limits_{i\in I}\alpha_i({\bf t}(i))^\top A_0{\bf t}(i)= {\bf e}^\top_k A_0 {\bftau}(j).\nonumber\eea
It is easy to see that these equalities can be rewritten as
$$ {\bf b}(k,j)=\left(\begin{array}{c}
A_m\bullet U\cr
m=0,...,n\end{array}\right)=\sum\limits_{i\in I}\alpha_i{\bf a}({\bf t}(i)) \mbox{ with } \alpha_i>0, \  {\bf t}(i)\in T,\  i \in I.$$
Thus we have  shown  that inclusions (\ref{cond-2}) hold true.

 \vspace{3mm}

$\Leftarrow)$   Now, having  supposed that  inclusions   (\ref{cond-2}) hold true, let us show that  the
 consistent system $\ca( {\bf x})\in \cop$ yields  the  uniform  LP  duality.

 Consider any ${\bf c}\in \mathbb R^n$ such that $Val({\bf P})>-\infty.$ It was stated  in section \ref{S-2} that  the  problem $({\bf P})$ is equivalent to the  problem
 (${\bf P}_*$) and
  there exists ${ \bf x}^*\in\mathbb R^n$ satisfying (\ref{W7}).

Notice that $Val({\bf P}_*)=Val({\bf P})$ and system (\ref{W6}) in  problem (${\bf P}_*$) can be rewritten as
\be (1, {\bf x}^{ \top}){\bf a}({\bf t})\geq 0 \ \forall {\bf t} \in \Omega,\  (1,{\bf x}^{\top}){\bf b}(k,j)  \geq 0 \ \forall k \in N_{*}(j),\ \forall j \in J.\label{W8}\ee

Taking into account the  inequalities in (\ref{W7}) (that can be considered as  a {\it generalized  Slater condition}),
 let us  show that system (\ref{W6})  yields  the  uniform LP  duality.
In fact, it follows from Theorem 1 in \cite{levin} that   under conditions (\ref{W7}), there exist vectors ${\bf t}(i)\in \Omega,\ i \in I, |I|\leq n,$ such that a Linear
Programming problem
$${\bf{LP}}:\ \ \min {\bf c}^\top{\bf x}  \ \mbox{ s.t. } (1,{\bf x}^\top){\bf a}({\bf t}(i))\geq 0 \ \forall i \in I,\  (1,{\bf x}^\top){\bf b}(k,j)
 \geq 0 \ \forall k \in N_{*}(j),\ \forall j \in J,$$
has the same optimal value as  the  problem (${\rm {\bf P}_*}$): $Val({\bf P}_*)=Val({\bf{LP}})>-\infty.$ The problem ({\bf{LP}}) is consistent since any ${\bf x}\in X$ is
feasible in this problem. Hence  the  problem ({\bf{LP}}) has an optimal solution.
Consequently,  there exist numbers and vector
$$\alpha_i,\ {\bf t}(i)\in \Omega,\ i \in I,   \ \lambda_k(j), \, k\in   N_{*}(j),\ j \in J,$$
such that $\alpha_i\geq 0,i \in I,$ $ \lambda_k(j)\geq 0, \, k\in N_{*}(j),\ j \in J,$ and
\be \sum\limits_{i\in I}\alpha_i{\bf a}({\bf t}(i))+ \sum\limits_{j\in J} \sum\limits_{k\in  N_*(j)}\lambda_k(j){\bf b}(k,j)=(-Val({\bf P}), c_m,\
m=1,...,n)^\top.\label{W12}\ee

From  (\ref{cond-2}), one can conclude  that for  all indices  $j\in J$, $k \in   N_{*}(j)$ and  any   $\lambda_k(j)>0$, the vector $\lambda_k(j){\bf b}(k,j)$ admits a
representation
\bea&\lambda_k(j){\bf b}(k,j)=\sum\limits_{i\in I(k,j)}\alpha_i(k,j){\bf a}({\bftau}(i,k,j))\nonumber\\
 &\mbox{ with } \alpha_i(k,j)>0,\, {\bftau}(i,k,j)\in T, \, i \in I(k,j), \, |I(k,j)|<\infty.\nonumber\eea
It follows from the representations  above  and  from  (\ref{W12}) that
\be \sum\limits_{i\in \bar I}\bar \alpha_i{\bf a}(\bar{\bf t}(i))=(-Val({\bf P}),c_m, m=1,...,n)^\top\ \mbox{ with some } \bar \alpha_i>0,\,  \bar{{\bf t}}(i)\in T, \, |\bar
I|<\infty. \label{W14}\ee
Denote
$U:=\sum\limits_{i\in \bar I}\bar \alpha_i\bar{\bf t}(i)(\bar{\bf t}(i))^\top.$
It is evident that $U\in \cp$ and relations (\ref{W14}) can be rewritten as  (\ref{W1}).
Hence we have shown that if   inclusions   (\ref{cond-2}) hold true, then  the  system $\ca({\bf x})\in \cop$ yields  the   uniform LP duality.
$\ \Box$

\vspace{3mm}

For $j\in J$ and $k\in N_*(j)$,  consider the following sets:
$$X(k,j):=\{{\bf x}\in X: {\bf e}^\top_k\ca({\bf x}){\bftau}(j)=0\}, \ T_{im}(k,j):=\{{\bf t} \in T:{\bf t}^\top\ca({\bf x}){\bf t}=0\ \forall{\bf x} \in X(k,j)\}.$$
  Denote $N(j)=N_*(j)\setminus M(j),$ $j\in J.$   Notice that by construction,
$$\emptyset\not =X(k,j)\subset  X\setminus\{{\bf x}^*\}, \  T_{im}\subset T_{im}(k,j)\  \forall   k\in N(j),\   \forall  j \in J;$$
$$  X(k,j)=X, \ \  T_{im}= T_{im}(k,j) \ \
 \forall k\in M(j),\   \forall j \in J.$$

  For  $j \in J$ and $k\in N_{*}(j)$, denote by ${\bf x}^*(k,j)$ a vector such that  ${\bf x}^*(k,j)\in X(k,j)$  and
\be  {\bf e}^\top_k\ca({\bf x}^*(k,j)){\bftau}(j)=0,\  {\bf t}^\top\ca({\bf x}^*(k,j)){\bf t}>0\ \forall {\bf t} \in T\setminus T_{im}(k,j).\label{W18}\ee
Notice that such vectors exist  for $j \in J, \ k\in N_{*}(j)$.

\begin{theorem}\label{P2}    A consistent linear  system $\ca({\bf x})\in \cop$
 with the corresponding sets $T_{im}$, ${\cal T}$, and other defined above,
 yields  the  uniform LP duality iff the following conditions hold:
\be{\bf I)}\qquad \qquad \quad  \   {\bf b}(k,j)\in {\rm cone}\{{\bf a}({\bf t}), {\bf t}\in T_{im}\}   \ \forall  k \in M(j),
\ \forall j \in J,\qquad \qquad \qquad \   \label{cond-a}\ee

\vspace{-8mm}

\be    {\bf II)}\qquad \quad \quad    {\bf b}(k,j)\in {\rm cone}\{{\bf a}({\bf t}), {\bf t}\in T_{im}(k,j)\}  \ \forall   k\in N(j),   \ \forall j \in J.\qquad \qquad \qquad
\label{cond-b}\ee
\end{theorem}

{\bf Proof.} It follows from Proposition \ref{P1} that to prove the theorem, it is enough to show that relations (\ref{cond-2}) are equivalent to relations (\ref{cond-a}) and
(\ref{cond-b}).

Since $T_{im}\subset T$ and $T_{im}(k,j)\subset T$ $\forall k\in N(j),$ $ j\in J$, it is evident  that the relations (\ref{cond-a}) and (\ref{cond-b}) imply  the   inclusions
(\ref{cond-2}).

 Suppose that   inclusions   (\ref{cond-2}) take a place. Hence   for any   $j\in J$ and  $k\in N_*(j),$ the equality
\be  {\bf b} (k,j)=\sum\limits_{i\in I}\alpha_i{\bf a}({\bf t}(i)) \mbox{ with some } \alpha_i=\alpha_i(k,j)>0,\,  {\bf t}(i)={\bf t}(i,k,j)\in T, \, i \in I,\label{W16}\ee
 holds true.
   Let's   multiply  the  right and left  parts  of  this   equality  by  $(1,({\bf x}^*(k,j))^\top)$.
As a result, we get
 \be \sum\limits_{i\in I}\alpha_i({\bf t}(i))^\top \ca({\bf x}^*(k,j)) {\bf t}(i)= {\bf e}^\top_k\ca({\bf x}^*(k,j)){\bftau}(j)=0 \ \mbox{ with some } \alpha_i>0,\, {\bf
 t}(i)\in T, i \in I.\label{06-12-4}\ee
 Here we took into account the equality in (\ref{W18}).   It follows from  (\ref{06-12-4})  and  the  inequalities in  (\ref{W18}) that
 in  (\ref{W16}), the  vectors ${\bf t}(i),i\in I,$  should  satisfy the conditions   ${\bf t}(i)\in T_{im}(k,j),$ $ i \in I.$ Consequently, we  have  shown  that inclusions (\ref{cond-2}) imply  the  inclusions
$${\bf b}(k,j)\in {\rm cone}\{{\bf a}({\bf t}), {\bf t}\in T_{im}(k,j)\},  k\in N_*(j)=M(j)\cup N(j), \ j \in J.$$
 Taking into account that $
T_{im}(k,j)=T_{im}$ $\forall  k \in M(j),$ $\forall  j\in J,$ we conclude that inclusions (\ref{cond-2}) imply
(\ref{cond-a}) and
(\ref{cond-b}). $\  \Box$

\vspace{3mm}

Notice that in (\ref{cond-a}) and (\ref{cond-b}) we have a {\it finite number} of inclusions.

 \section{Equivalent formulations of   the  condition I)}\label{s:4}

In this section, we   will present  several   equivalent formulations of   condition {\bf I}),  set forth  in the previous section, which is one of the   conditions that
guarantee the uniform  LP  duality of the copositive system. This  gives us the opportunity to  analyze this  condition from different points of view and   create a
theoretical basis  for comparing  our results with  others   known in the literature.

\begin{proposition}\label{P-new-1} Given a linear system $\ca({\bf x})\in \cop$   with the corresponding sets  $T_{im}, $ ${\cal T}, $ $M(j),\ j \in J$,
 and the vectors ${\bf a}({\bf t}), {\bf b}(k,j), k\in M(j),\ j \in J$ defined above, the following   statements  are equivalent:
\begin{enumerate}
\item[{j)}]  the  condition $\rm \bf{I})$  is satisfied;
\vspace{-2mm}

\item[{jj)}]  the    cones   ${\rm cone }\{{\bf b}(k,j),\, k \in M(j),\, j \in J\}$
  and   ${\rm cone }\{{\bf a}({\bf t}), {\bf t} \in T_{im}\} $   coincide;

  \vspace{-2mm}

\item[{jjj)}] the equality ${\rm span }\{{\bf b}(k,j),\, k \in M(j),\, j \in J\}={\rm cone }\{{\bf a}({\bf t}), {\bf t} \in T_{im}\}$  holds true.
\end{enumerate}
\end{proposition}

In  the  terminology from \cite{Duffin},   condition {\it jj)} means that the sets $\{{\bf b}(k,j),\, k \in M(j),\, j \in J\}$ and $\{{\bf a}({\bf t}), {\bf t} \in T_{im}\}$
are {\it positively
equivalent}.

\vspace{3mm}

{\bf Proof.} It is evident that  the  condition {\it jj)} implies  the  condition {\it j)}. Let us show that {\it j)} implies {\it jj)}. In fact, it follows from {\it j)}
that
${\rm cone }\{{\bf b}(k,j),\, k \in M(j),\, j \in J\}\subset {\rm cone }\{{\bf a}({\bf t}), {\bf t} \in T_{im}\}$. On the other hand, it follows from Proposition
\ref{Pr-07-12-1} that ${\rm cone }\{{\bf a}({\bf t}), {\bf t} \in T_{im}\}\subset {\rm cone }\{{\bf b}(k,j),\, k \in M(j),\, j \in J\}.$ Hence, ${\rm cone }\{{\bf b}(k,j),\, k
\in M(j),\, j \in J\} = {\rm cone }\{{\bf a}({\bf t}), {\bf t} \in T_{im}\} $, and  we have shown that {\it j)} implies that {\it jj)}. Thus the equivalence of {\it j)} and
{\it jj)} is proved.

 To prove the equivalence of   the  conditions {\it jj)} and   {\it jjj)}, it is enough to show  that
 \be {\rm cone }\{{\bf b}(k,j),\, k \in M(j),\, j \in J\}={\rm span }\{{\bf b}(k,j),\, k \in M(j),\, j \in J\}.\label{07-12-9}\ee
  From Proposition \ref{P0}, it follows  that $-{\bf b}(\bar k,\bar j)\in {\rm cone} \{{\bf b}(k,j),\, k \in M(j),\, j \in J\}$
	for all $\bar k\in M(\bar j)$ and $\bar j\in J$. This implies equality (\ref{07-12-9}) and hence  the  conditions {\it jj)} and {\it jjj)} are equivalent.
	 $\  \Box$

\vspace{3mm}

  Let the set $J$ be  partitioned  into  subsets $J(s),$ $ s \in S$, as in (\ref{W32}).
For $s \in S,$ let us number  the  indices in $J(s)$ as  follows:
$J(s)=\{i_1,i_2,...,i_{k(s)}\}, \ k(s)=|J(s)|$ and
denote
\bea& V(s):=\{(i_k,i_q), \ k=1,...,k(s),\, q=k,...,k(s)\},\ \ V_0:=\bigcup\limits_{s\in S}V(s).\label{V0}\eea
Notice that the  sets   introduced above  are finite.

For  a given vector  ${\bf z}=(z_0,z_1,...,z_n)^\top$, denote
\begin{equation}\cb({\bf z})=\sum\limits_{m=0}^nA_mz_m.\label{Bz}\end{equation}

\begin{proposition}\label{P4}   The condition $ {\rm \bf{I}})$  of Theorem \ref{P2}  (see (\ref{cond-a})) is equivalent to the following two conditions:

\begin{enumerate}

\item[{\bf A1)}] The set
$
L:={\rm cone}\{{\bf a}({\bf t}), {\bf t}\in T_{im}\}$
is a subspace;

\item[{\bf B1)}]  for any ${\bf z}\in \mathbb R^{n+1}$,   the  equalities
\be ({\bftau}(i))^\top\cb({\bf z}){\bftau}(j)=0  \ \forall    (i,j)\in V_0,\label{cond-bb-1}\ee
imply the equalities
\be {\bf e}^\top_k\cb({\bf z}){\bftau}(j)=0  \ \forall k\in M(j), \ \forall j\in J .\label{cond-bb-2}\ee
\end{enumerate}
\end{proposition}

{\bf Proof.} Suppose that   inclusions  (\ref{cond-a}) hold true. Then it follows from Proposition \ref{P-new-1} that $L$  is a subspace.

Now let us prove that   inclusions  (\ref{cond-a}) imply  the  condition {\bf{B1}}). First,
notice that it follows from Proposition \ref{P01} (see Appendix)   that   equalities (\ref{cond-bb-1}) are equivalent to the  equalities
\be {{\bf t}}^\top\cb({{\bf z}}){{\bf t}}=0  \ \forall {{\bf t}} \in T_{im}.\label{W30}\ee
There equalities and equalities (\ref{cond-bb-2}) can be written  as follows:
\be {\bf z}^\top {\bf a}({\bf t})=0  \ \forall {\bf t} \in T_{im},\label{W31}\ee
\be {\bf z}^\top {\bf b}(k,j)=0 \   \forall    k \in M(j),\    \forall   j \in J.\label{W33}\ee
Then it is evident that under conditions (\ref{cond-a}),  equalities (\ref{W31}) imply the equalities (\ref{W33}).
Thus we  have   shown  that  the  condition  {\bf{B1}}) follows from (\ref{cond-a}).

\vspace{2mm}

Now we will show that the conditions  {\bf{A1}}) and  {\bf{B1}}) imply   inclusions (\ref{cond-a}).

Notice that   under   the  condition  {\bf{A1}})   a vector  ${\bf z}$ satisfies
 (\ref{W31}) iff ${\bf z} \in L^\perp$    where $L^\perp$ is the orthogonal complement to $L$ in $\mathbb R^{n+1}.$    Hence it follows from   the  condition  {\bf{A1}}) and
 Proposition \ref{P01} that  {the} condition  {\bf{B1}}) can be reformulated as follows:
$${\bf e}^\top_k\cb({\bf z}){\bftau}(j)=0  \ \forall k\in M(j),  \ \forall j\in J, \ \forall {\bf z} \in L^\bot,$$
or equivalently
\be {\bf z}^\top {\bf b}(k,j)=0  \ \forall k\in M(j), \ \forall j\in J, \ \forall {\bf z} \in L^\bot,\label{W34}\ee
  Given  $k\in M(j)$ and $j\in J$, the vector ${\bf b}(k,j)$ admits the representation
$${\bf b}(k,j)={\bf b}^1(k,j)+{\bf b}^2(k,j) \ \mbox{ with } {\bf b}^1(k,j)\in L \mbox{ and } {\bf b}^{2}(k,j)\in L^\bot.$$
It follows from (\ref{W34}) that ${\bf b}^2(k,j)=0$ and ${\bf b}(k,j)={\bf b}^1(k,j)\in L$  in the representation above.
Hence  the  conditions  {\bf{A1}}) and  {\bf{B1}}) imply   the inclusions (\ref{cond-a}).  $\  \Box$
\begin{remark} It is easy to see that   under  the  condition {\bf{A1}}), the condition  {\bf{B1}}) can be reformulated as  {follows}:\  the equalities (\ref{cond-bb-2}) hold
true for any ${\bf z}\in L^\perp.$
\end{remark}
 \begin{remark} Let us introduce $n+1$ vectors
${\bf a}(i,j):=((\bftau(i))^\top A_m\bftau(j), \, m=0,...,n)^\top, \, (i,j)\in V_0,$
and consider matrices with finite dimensions
$$\mathbb A:=({\bf a}(i,j), \, (i,j)\in V_0),\ \mathbb B:=({\bf b}(k,j),\ k \in M(j),\, j \in J).$$
Then the condition {\bf B1}) can be formulated as
$\ {\rm rank}\mathbb A={\rm rank}(\mathbb A,\, \mathbb B).$
\end{remark}

\begin{proposition}\label{P5}  The condition  {\bf{A1}}) is equivalent to the following  {one:}

\vspace{2mm}

{\bf A2)} There exists a matrix $U^*\in \cp$ in the form
\bea& U^*=\sum\limits_{i \in I}\alpha_i{\bf t}(i)({\bf t}(i))^\top+\frac{1}{4}\sum\limits_{(l,q)
\in V_0}({\bftau}(l)+{\bftau}(q))({\bftau}(l)+{\bftau}(q))^\top \label{U*}\\&\mbox{ with } \alpha_i>0, \,{\bf t}(i)\in T_{im}, \,
 i \in I,\  |I|<\infty, \nonumber\eea
such that
\be A_m\bullet U^*=0 \ \forall m=0,1,...,n.\label{W40}\ee
\end{proposition}
{\bf Proof.} Suppose that  {the} condition  {\bf{A1}})  holds true.  Hence it follows from Proposition \ref{P02} (see Appendix)  that there exist numbers and vectors
 $\bar{\alpha}_i>0, \ {\bf t}(i)\in T_{im}, i \in I, \ |I|<\infty,$ such that
 $$ \sum\limits_{i \in I}\bar \alpha_i{\bf a}({\bf t}(i))=0,\ {\rm rank} ({\bf a}({\bf t}(i)), i \in I)=r_{im},$$
 where $r_{im}:={\rm rank} ({\bf a}({\bf t}), {\bf t} \in T_{im}).$    Notice that $\frac{1}{2}(\bftau(l)+\bftau(q))\in T_{im}$ for all $(l,q)\in V_0.$ Then, according to
 Proposition \ref{P03} (see the Appendix), there exist numbers $\alpha_i>0, i \in I,$ such that
\be \sum\limits_{i \in I} \alpha_i{\bf a}({\bf t}(i))+\frac{1}{4}\sum\limits_{(l,q)\in V_0} {\bf a} ( {\bftau}(l)+{\bftau}(q))=0.\label{W41}\ee

Let $U^*$ be  a matrix in the form (\ref{U*}). Then (\ref{W41}) can be rewritten in the form (\ref{W40}). Thus we have proved that  {the} condition {\bf{A1}}) implies the
condition  {\bf{A2}}).

Now suppose that the condition  {\bf{A2}}) holds true   and rewrite equalities  (\ref{W40}) in the form (\ref{W41}).
  It follows from Proposition \ref{P04}  (see Appendix)   that
$${\rm rank} ({\bf a}({\bf t}(i)), i \in I, \ {\bf a}({\bftau}(l)+{\bftau}(q)),\  (l,q)\in V_0)=r_{im}. $$
Taking into account this equality,  equality (\ref{W41}), and Proposition \ref{P02}, we conclude
that $L:={\rm cone}\{{\bf a}({\bf t}), {\bf t} \in T_{im}\} $ is a subspace. Thus we have shown that  {the}
 condition {\bf{A2}}) implies the condition  {\bf{A1}}). $\ \Box$

\begin{remark}  The condition  {\bf{A2}}) can be reformulated as follows: there
exists a matrix $U^*\in \cp$ in the form
$$U^*=\sum\limits_{i \in I}\alpha_i{\bf t}(i)({\bf t}(i))^\top\mbox{ with } \alpha_i>0, {\bf t}(i)\in T_{im},  i \in I,$$
such that ${\rm rank}({\bf a}({\bf t}(i)), \,i \in I)=r_{im}$ and equalities
(\ref{W40}) hold true.
\end{remark}

\vspace{3mm}

To formulate the next   propositions and   lemmas, we need the following notation and definitions (see \cite{Pataki2017}).  Given a matrix  $Y\in {\cal S}^p, $ denote
$${\rm dir}(Y, \cop): =\{ D\in {\cal S}^p:  \ Y + \varepsilon D \in \cop \mbox{ for some } \varepsilon  >0 \},$$
$${\rm ldir}(Y, \cop): = {\rm dir}(Y, \cop) \cap  -{\rm dir}(Y, \cop),$$
$${\rm tan}(Y, \cop): = {\rm cl}\,  ( {\rm dir}(Y, \cop) ) \cap  -{\rm cl}\, ( {\rm dir}(Y, \cop)).$$

 For     matrices  $Y\in \cop $  and $U\in \cp$,  we say that $U$ is {\it strictly complementary} to $ Y $ if $U\in {\rm relint} (\cp\cap Y^\perp)$.

 \begin{definition}\label{def2} A matrix $Y$ is called a {\it maximum slack} in the system $\ca({\bf x})\in \cop$ if $Y\in  {\rm relint}\, {\cal D}$, where ${\cal
 D}:=\{D=\ca({\bf x}), \,{\bf x} \in X\}.$
 \end{definition}
Denote
$$ {\cal R}(\cb):=\{D=\cb({\bf z}), \, {\bf z} \in \mathbb R^{n+1}\},\ {\cal N}(\cb^*):=\{U\in {\cal S}^p: A_m\bullet U=0 \ \forall m=0,1,...,n\}$$
 where $\cb({\bf z})$ is defined in (\ref{Bz}).

\begin{lemma}\label{L1} A  {matrix}  $Y^*$  is a maximum slack  in  {the} system $\ca({\bf x})\in \cop$   iff there exists ${\bf y}^*\in X$ such that $Y^*=\ca({\bf y}^*)$ and
\be {\bf t}^\top \ca({\bf y}^*){\bf t}>0\ \forall {\bf t} \in T\setminus T_{im},\ {\bf e}^\top_k \ca( {\bf y}^*){\bftau}(j)>0\ \forall k \in P\setminus M(j), j \in
J.\label{ymin}\ee
\end{lemma}
{\bf Proof.}
Suppose that relations (\ref{ymin}) hold true   for some ${\bf y}^*\in X$. Let us
  show that $Y^*:=\ca({\bf y}^*)\in {\rm relint}\,{\cal D}. $  Here, as above, ${\cal D}:=\{D=\ca({\bf x}), \,{\bf x} \in X\}.$
Since  the set ${\cal D}$ is convex{, we can state that the following equivalence takes a place:}
\bea& \ca({\bf y}^*)\in {\rm relint}\,{\cal D}\ \ \Longleftrightarrow\  \forall{\bf x} \in X\ \exists\, \ep>0
 \mbox{ such that } D(\ep,{\bf x})\in \cop, \label{equiv}\\
& D(\ep,{\bf x}):=(1+\ep)\ca({\bf y}^*)-\ep \ca({\bf x}).  \nonumber\eea
Notice that due to Theorem 1 from \cite{KT-new}, we have
\bea&  D(\ep,{\bf x})\in \cop \ \Longleftrightarrow\label{12-12-22}\\
&{\bf  t}^\top D(\ep,{\bf x}) {\bf t}\geq 0 \; \forall\, {\bf t} \in \Omega; \
 D(\ep,{\bf x}) \bftau(j)\geq 0\ \forall j \in J.\label{07-12-11}\eea
  where the set $\Omega$ is defined in (\ref{Omega}).
Denote $$\xi_1:=\min\limits_{{\bf t}\in \Omega}t^\top \ca({\bf y}^*)t,\
\xi_2:= \min\{{\bf e}^\top_k \ca({\bf y}^*)\bftau(j), \, k \in P\setminus M(j), \, j\in J\}.$$ It follows from (\ref{ymin}) that $\xi_1>0$ and $\xi_2>0.$

For a fixed ${\bf x}\in X,$ set
$$\eta_1({\bf x}):=\max\limits_{{\bf t}\in \Omega} {\bf t}^\top \ca({\bf x}) {\bf t},\ \eta_2({\bf x})=\max\{{\bf e}^\top_k \ca({\bf x})\bftau(j), \, k \in P\setminus M(j), \,
j\in J\},$$
 and  calculate
\bea & {\bf t}^\top D(\ep, {\bf x}){\bf t}=(1+\ep){\bf t}^\top \ca({\bf y}^*){\bf t}-
\ep {\bf t}^\top\ca({\bf x}){\bf t}\geq (1+\ep)\xi_1-\ep\eta_1({\bf x})\ \forall t \in \Omega,\nonumber\\
& {\bf e}^\top_k D(\ep, {\bf x}){\bftau}(j)=(1+\ep){\bf e}^\top_k \ca({\bf y}^*){\bftau}(j)-\ep {\bf e}^\top_k\ca({\bf x}){\bftau}(j)\geq
 (1+\ep)\xi_2-\ep\eta_2({\bf x})\label{07-12-12}\\ & \qquad \qquad\qquad \forall \,
k\in P\setminus M(j),\ \forall j \in J.\nonumber\eea
For  $s=1,2$, set
$\ep_s({\bf x})=\infty  \mbox{ if } \xi_s-\eta_s({\bf x})\geq 0;\ \ep_s({\bf x})=\xi_s/(\eta_s({\bf x})-\xi_s)>0  \mbox{ if } \xi_s-\eta_s({\bf x})< 0;$
$\ep_0({\bf x}):=\min\{\ep_1,\ep_2\}>0.$
It follows from (\ref{07-12-12}) that for any ${\bf x}\in X,$ there exists $\ep_0({\bf x})>0$ such that inequalities (\ref{07-12-11}) hold true. Hence we
have shown that inequalities (\ref{ymin}) imply the condition $\ca({\bf y}^*)\in {\rm relint}\,{\cal D}.$

Now suppose that {$Y^*$  is a maximum slack  in  {the} system $\ca({\bf x})\in \cop$}.  {Then} there exists ${\bf y}^*\in X$
 such that $Y^*=\ca({\bf y}^*)$,  {$\ca({\bf y}^*)\in {\rm relint}\, {\cal D},$} and due to (\ref{equiv}) we have
\be \forall \,{\bf x} \in X  \ \exists \ \ep>0\  \mbox{ such that }  {D(\ep, {\bf x})}\in \cop .\label{11*}\ee

Suppose that  this  vector  ${\bf y}^*$   does not satisfy (\ref{ymin}). Consequently one of the following situations should have place:

\begin{enumerate}

\item[\textbf{a})] there exists $\bar{\bf t} \in T\setminus T_{im}$ such that $\bar{\bf t}^\top \ca({\bf y}^*)\bar{{\bf t}} =0;$

\item[\textbf{b})] ${\bf t}^\top \ca( {\bf y}^*){\bf t} $ for all $  {\bf t} \in T\setminus T_{im}$ and there exist $j_0\in J$ and $k_0\in P\setminus M(j_0)$ such that
$ {\bf e }^\top_{k_0}\ca({\bf y}^*){\bftau}(j_0)=0.$
\end{enumerate}
Let ${\bf x}^*$  be  a vector in $X$ satisfying (\ref{W7}). We will show that for  ${\bf x}={\bf x}^*$, inclusion   (\ref{12-12-22}) does not hold true with any $\ep>0.$


In fact, in  {the} situation \textbf{a}), we have that $\bar{\bf t}^\top \ca({\bf y}^*)\bar{{\bf t}} =0$ and
$\bar{\bf t}^\top \ca({\bf x}^*)\bar{{\bf t}} >0$, and hence
$$\bar{\bf t}^\top ( (1+\ep)\ca({\bf y}^*)-\ep \ca({\bf x}^*))\bar{\bf t} <0 \mbox{ with any } \ep>0.$$
 It follows from these relations and  the  condition $\bar {\bf t}\in T$  that
\be (1+\ep)\ca({\bf y}^*)-\ep \ca({\bf x}^*)= {D(\ep, {\bf x}^*)}\not \in \cop \ \forall \ep>0.\label{07-12-18}\ee

In  the  situation \textbf{b}), for $\theta \geq 0$ and
${\bf t}(\theta):=( {\bftau}  (j_0)+\theta{\bf e}_{k_0})\in \mathbb R^p_+,$ let us calculate
$${\bf t}^\top(\theta) {D(\ep, {\bf x}^*)}{\bf t}(\theta)
=[2\theta{\bf e}^\top_{k_0} {D(\ep, {\bf x}^*)}{\bftau}(j_0)+
\theta^2{\bf e}^\top_{k_0} {D(\ep, {\bf x}^*)}{\bf e}_{k_0}]$$
$$=\theta[-2\ep{\bf e}^\top_{k_0} \ca({\bf x}^*){\bftau}(j_0)+\theta{\bf e}^\top_{k_0}((1+\ep)\ca({\bf y}^*)-\ep \ca({\bf x}^*)){\bf e}_{k_0}] $$
\be  \leq  \theta[-2\ep{\bf e}^\top_{k_0} \ca({\bf x}^*){\bftau}(j_0)+\theta(1+\ep){\bf e}^\top_{k_0}\ca({\bf y}^*){\bf e}_{k_0}].\label{07-12-16}\ee
For any $\ep>0$, let us set
$\theta(\ep) =\ep{\bf e}^\top_{k_0} \ca({\bf x}^*)){\bftau}(j_0)/(1+\ep){\bf e}^\top_{k_0}\ca({\bf y}^*){\bf e}_{k_0}>0$ if ${\bf e}^\top_{k_0}\ca({\bf y}^*){\bf e}_{k_0}>0$
and $\theta(\ep)=1$ if ${\bf e}^\top_{k_0}\ca({\bf y}^*){\bf e}_{k_0}=0$.
Then taking into account (\ref{07-12-16}), for ${\bf t}(\theta(\ep))\in \mathbb R^p_+$ we obtain
$${\bf t}^\top(\theta(\ep)) {D(\ep, {\bf x}^*)}{\bf t}(\theta(\ep))\leq -\ep \theta(\ep){\bf e}^\top_{k_0} \ca({\bf x}^*){\bftau}(j_0)<0 \ \forall \ep>0.$$
This implies that relations (\ref{07-12-18}) take place.  But these relations contradict (\ref{11*}).
Hence we  have shown that (\ref{11*}) implies  (\ref{ymin}).$\  \Box$

\begin{proposition}\label{P6}  The condition   {\bf A1)} is equivalent to the following  one:

\vspace{2mm}

{\bf A3)} For  a maximum slack $Y^*$  in  {the} system $\ca({\bf x})\in \cop,$
there is $U^* \in {\cal N}(\cb^*) \cap \cp $ strictly complementary to $Y^*$.
 \end{proposition}
 {\bf Proof.} Suppose that  the  condition  {\bf A3}) holds true. Since
the set $\cp\cap (Y^*)^\perp$ is convex, then
 \bea& U^*\in {\rm relint} (\cp\cap (Y^*)^\perp)\, \ \Longleftrightarrow\, \nonumber\\
  & \forall U\in \cp\cap (Y^*)^\perp\  \ \exists\, \varepsilon>0\mbox{ such that } (1+\varepsilon)U^*-\varepsilon U\in \cp\cap (Y^*)^\perp.\label{W52}\eea
 It follows from Lemma \ref{L1} that  $U\in \cp\cap (Y^*)^\perp$ iff $U$ admits a representation
 \be U=\sum\limits_{i \in \bar I}\bar \alpha_i\bar{\bf t}(i)(\bar{\bf t}(i))^\top \ \mbox{with some } \bar \alpha_i>0,\
 \bar{\bf t}(i)\in T_{im}, \ i \in \bar I.\label{W73}\ee
Hence $U^*$ takes the form  \be U^*=\sum\limits_{i \in I}\alpha_i {\bf t}(i)({\bf t}(i))^\top \mbox{ with some } \alpha_i>0,\  {\bf t}(i)\in T_{im}, \ i \in  I.
\label{W71}\ee
Since $U^* \in {\cal N}(\cb^*)$ then  the equalities (\ref{W40}) hold true.

For any $ {\bf t}\in T_{im}$, consider  the  matrix $U={\bf t}{\bf t}^\top\in \cp\cap (Y^*)^\perp.$  It follows from (\ref{W52}),  (\ref{W73})
that there exist  $\varepsilon>0$ and
$\bar \alpha_i>0,$ $\bar{\bf t}(i)\in T_{im}$, $i \in \bar I,$
such that
$$(1+\varepsilon)U^*-\varepsilon {\bf t}{\bf t}^\top=\sum\limits_{i \in \bar I}\bar \alpha_i\bar{\bf t}(i)(\bar{\bf t}(i))^\top \
\Longleftrightarrow \
(1+\varepsilon)U^*=\varepsilon {\bf t}{\bf t}^\top+\sum\limits_{i \in \bar I}\bar \alpha_i\bar{\bf t}(i)(\bar{\bf t}(i))^\top.$$
Taking into account the latter  equality    and  (\ref{W40}),    we obtain   the equalities
$$A_m\bullet (\varepsilon {\bf t}{\bf t}^\top+\sum\limits_{i \in \bar I}\bar \alpha_i\bar{\bf t}(i)(\bar{\bf t}(i))^\top)=0\ \forall m=0,1,...,n,$$
  that can be rewritten in the form
$$\varepsilon  {\bf a}({\bf t})+\sum\limits_{i \in \bar I}\bar \alpha_i{\bf a}(\bar{\bf t}(i))= 0 \mbox{ where } \varepsilon>0, \,\bar \alpha_i>0, \, \bar{\bf t}(i)\in T_{im}
\,\forall i \in \bar I.$$
It follows from these relations that $-{\bf a}({\bf t})\in L$ for any ${\bf a}({\bf t})\in L$ and hence, $L$ is a subspace.

\vspace{3mm}

Now suppose that the condition \textbf{A1}) holds true. It follows from
 Proposition \ref{P5} that the condition {\bf A1}) is equivalent to the condition {\bf A2}). According to this condition,
there exists a matrix $U^*\in \cp$ in the form (\ref{U*}) satisfying equalities (\ref{W40}).
By construction, $U^*\in \cp\cap (Y^*)^\perp$ and $U^* \in {\cal N}(\cb^*) \cap \cp $.
Let us show that relations (\ref{W52}) hold true.

Consider any  {matrix} $U\in \cp\cap (Y^*)^\perp$.    It follows from Lemma \ref{L1} that   this matrix  admits  representation (\ref{W73}).
For a fixed $i \in \bar I$, consider   the corresponding   $\bar{\bf t}(i)\in T_{im}$.  Then it follows from Proposition \ref{P05}  (see Appendix)  that the matrix
$\bar{\bf t}(i)(\bar{\bf t}(i))^\top$   can be presented in the form
$$\bar{\bf t}(i)(\bar{\bf t}(i))^\top=\sum\limits_{(l,q) \in V_0}\beta_{l,q}(i)({\bftau}(l)+{\bftau}(q))({\bftau}(l)+{\bftau}(q))^\top.$$
Consequently,
\begin{equation*}\begin{split}(1+\varepsilon)U^*-\varepsilon U=&(1+\varepsilon)\Bigl(\sum\limits_{i \in I}\alpha_i{\bf t}(i)({\bf t}(i))^\top+\frac{1}{4}\sum\limits_{(l,q) \in
V_0}({\bftau}(l)+{\bftau}(q))({\bftau}(l)+
{\bftau}(q))^\top\Bigl)\\
& {-}\varepsilon\sum\limits_{i\in \bar I}\bar \alpha_i\Bigl(\sum\limits_{(l,q) \in V_0}\beta_{l,q}(i)({\bftau}(l)+{\bftau}(q))( {\bftau}(l)+{\bftau}(q))^\top\Bigl)\\
& {=}(1+\varepsilon)\sum\limits_{i \in I}\alpha_i{\bf t}(i)({\bf t}(i))^\top+\sum\limits_{(l,q) \in V_0}\bar
\beta_{l,q}({\bftau}(l)+{\bftau}(q))({\bftau}(l)+{\bftau}(q))^\top,
\end{split}\end{equation*}
where $\alpha_i>0, \, {\bf t}(i)\in T_{im} \, \forall i \in I, \ 0.5({\bftau}(l)+\ {\bftau}(q))\in T_{im}, \, \forall (l,q)\in V_0,$
 and  {where} for  a sufficiently small $\varepsilon>0$,  {it holds} $\bar \beta_{l,q}:=(1+\varepsilon)/4-\varepsilon\sum\limits_{i\in \bar I}\bar \alpha_i\beta_{l,q}(i)>0 \
 \forall (l,q)\in V_0$.

Then, evidently,  $(1+\varepsilon)U^*-\varepsilon U\in \cp\cap (Y^*)^\perp $ for any $U \in \cp\cap (Y^*)^\perp$
and for a some sufficiently small $\varepsilon>0$, and, consequently $U^*\in {\rm relint} (\cp\cap (Y^*)^\perp)$.
  $\ \Box$

 \begin{proposition}\label{P7}  {The c}ondition {\bf B1}) is equivalent to the following  {one:}

 \vspace{2mm}

 {\bf B2)}  For  a maximum slack $Y^*$  in  {the} system $\ca({\bf x})\in \cop$, the set
 $$ {\cal R}(\cb) \cap ({\rm tan}(Y^*, \cop) \setminus {\rm ldir}(Y^*, \cop))$$
 is empty.
\end{proposition}
{\bf Proof.}
In \cite{Dick} (see Theorems 6, 13), for $A\in \cop$, it is shown that
\bea {\rm dir}(A,\cop)=&&\{B\in {\cal S}^p: {\bf t}^\top B{\bf t}\geq 0 \ \forall {\bf t} \in {\cal V}^A;\nonumber\\
&&{\bf e}^\top_kB{\bf t}\geq 0\ \forall {\bf t} \in {\cal V}^A\cap {\cal V}^B,\ \forall k\in \{s \in P:{\bf e}^\top_sA{\bf t}=0\}\},\label{dirA}\eea
\bea {\rm tan}(A,\cop)=&&\{B\in {\cal S}^p: {\bf t}^\top B{\bf t}=0\ \forall {\bf t} \in {\cal V}^A\}=\label{tanA}\\
&&\{B\in {\cal S}^p: {\bf t}^\top B{\bftau}=0\ \forall \{{\bf t},{\bftau}\} \subset  {\cal V}_{min}^A \ \mbox{s.t.} \ {\bf t}^\top A {\bftau} =0\},\nonumber\eea
where ${\cal V}^A:=\{{\bf t} \in T: {\bf t}^\top A {\bf t}=0\}$ is the set of zeros of $A$ and ${\cal V}_{min}^A \subset {\cal V}^A $ is the set of minimal zeros of $A$. For
definitions see \cite{Dick}.

Let us   present    the  condition {\bf B2}) in another form. It follows from Lemma \ref{L1} that if $Y^*$  is  a maximum slack  in  the  system $\ca({\bf x})\in \cop$, then
there exists  {a vector} ${\bf y}^*\in X$ such that   {conditions} (\ref{ymin}) hold true. Taking into account {the relations} (\ref{ymin}), (\ref{dirA}),  and (\ref{tanA}),
it is easy to see that
$${\rm tan}(\ca({\bf y}^*),\cop)\setminus {\rm ldir}(\ca({\bf y}^*),\cop)=\{B\in {\cal S}^p:({\bftau}(i))^\top B{\bftau}(j)=0  \,\forall (i,j)\in V_{0};$$
$$\exists \, \bar{\bf t} \in T_{im}\cap {\cal V}^{B} \mbox{ and } \exists \, \bar k\in P\mbox{ such that }{\bf e}^\top_{\bar k}\ca({\bf y}^*)\bar{\bf t}=0,\ {\bf e}^\top_{\bar
k}B\bar{\bf t}\not =0\}.$$
Consequently,  the  condition  {\bf B2)}
 can be reformulated as follows:

\vspace{2mm}

 {\bf B2*): {\it for any} ${\bf z} \in \mathbb R^{n+1},$} {\it the  equalities}
 $({\bftau}(i))^\top \cb({\bf z}){\bftau}(j)=0\ \forall (i,j)\in V_{0}$ {\it imply the equalities}
\be {\bf e}^\top_{ k}\cb({\bf z}){\bf t} =0\ \forall \, k \in \{q\in P: {\bf e}^\top_{ q}\ca({\bf y}^*){\bf t} =0\}
\ \forall  {\bf t} \in T_{im}\cap {\cal V}^{\cb({\bf z})}=T_{im}.  \label{W74}\ee

Suppose that that  the  condition  {\bf B2*)}  holds  true. For $j \in J$,
consider the corresponding  vector  ${\bftau}(j)\in T_{im}$. By construction (see (\ref{ymin})), we have
$\{q\in P: {\bf e}^\top_{ q}\ca({\bf y}^*){\bftau}(j) =0\} =M(j).$ Consequently, it follows from conditions (\ref{W74}) that
$ {\bf e}^\top_{ k}\cb({\bf z}){\bftau}(j)=0$ for all $k \in M(j).$ Hence we have shown  that  the  condition  {\bf B2*)}  implies the  condition
{\bf{B1})}.

Now suppose that  {the} condition  {\bf{B1}}) holds true.
Consider any ${\bf t}\in T_{im}$. It follows from (\ref{W32})
 that ${\bf t} \in T_{im}(s)$ with some $s \in S$ and consequently, ${\bf t}$ admits a representation
\be {\bf t}=\sum\limits_{j\in \Delta J(s)}\alpha_j{\bftau}(j),\ \alpha_j>0, \, j\in \Delta J(s)\subset J(s); \ \sum\limits_{j\in\Delta J(s)}\alpha_j=1.\label{3*}\ee

Suppose that  $k\in \{q\in P: {\bf e}^\top_{ q}\ca({\bf y}^*){\bf t} =0\}$. Hence
$$ 0={\bf e}^\top_{k}\ca({\bf y}^*){\bf t}=\sum\limits_{j\in\Delta J(s)}\alpha_j{\bf e}^\top_{k}\ca({\bf y}^*){\bftau}(j).$$
Taking into account this equality and  the  inequalities $\alpha_j>0, \, j\in \Delta J(s)$, and ${\bf e}^\top_q\ca({\bf y}^*){\bftau}(j)>0$ $ \forall q \in P\setminus M(j)$
$\forall j \in J(s),$ we conclude that
\be k \in M(j) \ \forall j \in \Delta J(s).\label{4*}\ee

Now, for the same ${\bf t} \in T_{im}(s)$  and any ${\bf z}\in \mathbb R^{n+1}$ satisfying (\ref{cond-bb-1})
calculate ${\bf e}^\top_{k}\cb({\bf z}){\bf t}$ taking into account  conditions (\ref{cond-bb-2}), (\ref{3*}), and (\ref{4*}):
$${\bf e}^\top_{k}\cb({\bf z}){\bf t}=\sum\limits_{j\in\Delta J(s)}\alpha_j{\bf e}^\top_{k}\cb({\bf z}){\bftau}(j) =0.$$
Thus we have  shown  that  the  condition  {\bf{B1}}) implies  {\bf{B2*}}).    $\ \Box$

\vspace{5mm}

The above considerations can be formulated as follows.
\begin{lemma} \label{L-nec} For any $k\in \{1,2,3\}$ and any $m\in \{1,2\}$,  {the} conditions i) and ii) below
are equivalent to each other, and  are  necessary for  the   consistent  system $\ca({\bf x})\in \cop$  to yield  the   uniform LP duality  property:

{ i)}   the condition $ {\rm \bf{I}})$ holds true;

{ ii)}   the conditions   {\bf {A}$k$)}  and   {\bf{B}$m$)} hold true.
\end{lemma}

  The condition $ii)$   with   $k=3$ and $m=2$  (i.e., the conditions {\bf A3}) and {\bf B2})) is a necessary condition for  the  linear conic system  to yield  the   uniform
  LP  duality formulated and proved in \cite{Pataki2017} and applied to  the copositive cone $\cop.$

It was shown in \cite{Pataki2017} that   when ${\cal K}$ is a nice cone, the conditions {\bf A3}) and {\bf B2}) are necessary and sufficient for   the 	linear conic
consistent  system  $\ca({\bf x})\in {\cal K}$ to  yield   the  uniform LP duality.

 In general,  the conditions formulated in Lemma \ref{L-nec}  are  only necessary, but not sufficient for   the  system $\ca({\bf x})\in \cop$  to yield  the   uniform LP
 duality. It is illustrated by a simple  example presented in the next section.

 \section{Examples}\label{s:5}

 In this section, we consider several examples that illustrate our results and help us to compare our results with ones known in   the  literature.

\vspace{2mm}

{\bf Example 1.}
Consider  the  system $\ca({\bf x})\in \cop$  with the following data:
{\small \begin{equation}  n=2,\; p=3, \   A_0=\left(\begin{array}{ccc}
a& 0 &0\\
0&0&0\\
0&0&0\end{array}\right),
A_1=\left(\begin{array}{ccc}
0& 0 &0\\
0&-1&0\\
0&0&0\end{array}\right),
A_2=\left(\begin{array}{ccc}
-1& 0 &0\\
0&0&-1\\
0&-1&0\end{array}\right), \ a>0.\label{data21}\end{equation}}

 For ${\bf t}^*=(0,0,1)^\top$, we have $({\bf t}^*)^\top{\cal A}({\bf x}){\bf t}^*=0,\; {\bf e}_1^\top{\cal A}({\bf x}){\bf t}^*=0$, $ {\bf e}_{3}^\top{\cal A}({\bf x}){\bf
 t}^*=0$ for all ${\bf x}\in \mathbb R^2.$

 It is easy to check that for the vector $ {\bf x}^*=(-1,-1)^\top$  and the data  in  (\ref{data21}) we  have
 ${\bf t}^\top{\cal A}({\bf x}^*){\bf t}>0$ for all ${\bf t}\in \mathbb R^3_+\setminus\{{\bf t^*}\}$ and ${\bf e}_2^{\top}{\cal A}({\bf x}^*){\bf t^*}>0.$ Hence, for the
 system
 under consideration,  $\ca({\bf x}^*)$ is a maximum slack, $T_{im}=\{ {\bftau}(1)={\bf t}^*\},$ $M(1)=\{1,3\}$, $J=\{1\},$
${\bf a}({\bf t}^*)=(0,\, 0,\,0)^\top, $ $ {\bf b}(1,1)={\bf b}(3,1)=(0,\, 0,\,0)^\top, $ ${\bf b}(2,1)=(0,\, -1,\,0)^\top.$
 Hence $L:={\rm cone}\{{\bf a}({\bf t}), {\bf t} \in T_{im}\}=\{  {\bf a}({\bf t}^*)\}$ and ${\bf b}(k,j)\in L$ for all $k \in M(j)$ and all $j \in J.$

 Thus we see that  for this system,  the  condition  {\bf I}) is satisfied,
	and it follows from Lemma \ref{L-nec} that  the conditions  {\bf A$k$}) for $k=1,2,3$ and  the  conditions   {\bf B$m$})
	for $m=1,2$ are satisfied as well.  (The fulfillment  of  the conditions  {\bf A$k$}) for $k=1,2,3$ and  the  conditions   {\bf B$m$})
	for $m=1,2$ can be checked  directly.)

  But the  system under consideration does not  yield  the   uniform  {LP} duality.
  In fact, it was shown in \cite{KTD} that for  the  primal problem ({\bf P}) with  the cost vector  $c^{\top}=(0,\ -1)$ and the corresponding dual problem ({\bf D}), there is
  a  positive  duality gap: $Val({\rm{\bf P}})-Val({\rm{\bf D}})=a>0.$

 The reason for   not complying with the uniform duality   is that for the  system with  the  data (\ref{data21}),
 the condition {\bf II}) is not satisfied.   Indeed,  for the vector $\bar{\bf x}=(-1,\, 0)^\top$,
we have ${\bf t}^\top \ca(\bar{\bf x}){\bf t}>0$ $ \forall {\bf t} \in T\setminus \{{\bf t}^*\}$ and
${\bf e}^\top _k\ca(\bar{ \bf x}){\bf t}^*=0$ for $k=1,2,3.$
 Hence for  the   system under consideration, we have $N(1)=\{2\},$  $T_{im}(k=2,j=1)=T_{im}=\{{\bf t}^*\}$  and  the  condition {\bf{II}}), for $j=1\in J$ and $k=2\in N(j)$,
 takes the form
 $${\bf b}(2,1)=(0,\, -1,\,0)^\top\in {\rm cone }\{{\bf a}({\bf t}), {\bf t} \in T_{im}(k,j)\}=\{{\bf a}({\bf t}^*)\}=\{(0,\, 0,\,0)^\top\}.$$
  It is evident that this condition  does not hold true.

 \vspace{3mm}

 This example shows that   the  condition {\bf{II}})  is essential and  can not be omitted.
  The example  also shows that for  the  cone $\cop$,  the conditions formulated in \cite{Pataki2017}
are not sufficient   unlike  the  case  with the   cone  of   positive semi-definite matrices ${\cal S}^p_+$   for which these conditions are necessary and sufficient.

\vspace{3mm}

{\bf Example 2.}  Let us  now  consider  an   example  where the  condition  {\bf I})  is violated.  Consider  the  system $\ca({\bf x})\in \cop$ with the following data:
\be n=1, \ p=3, \ A_0=
\left(\begin{array}{ccc}
a& 0 &-a\cr
0 & 0& 0\cr
-a& 0 &a\end{array}\right),  \
A_1=
\left(\begin{array}{ccc}
1& -1 & 2\cr
-1 & 0& 1\cr
2& 1 &-5 \end{array}\right),\  a>0.\label{80}\ee
This system admits a unique feasible solution ${\bf x}=x_1=0.$
Hence $X=\{{\bf 0}\}$ and it is easy to check that $T_{im}=\{{\bf t} \in T:t_1=t_3\}$,
where $T=\{{\bf t} \in \mathbb R^3_+:t_1+t_2+t_3=1\}.$  The vertices of the set ${\rm conv} T_{im}$
 are ${\bftau}(1)=0.5(1,\, 0,\, 1)^\top,$ ${\bftau}(2)=(0,\, 1,\, 0)^\top,$ and the sets $M(j)$, $N(j)=N_*(j)\setminus M(j),$
defined in (\ref{Mj}) and (\ref{W5}) take   the   form $M(j)=\{1,2,3\}$, $N(j)=\emptyset$ for $j \in J=\{1,2\}.$

It is easy to see that ${\bf t}^\top A_0{\bf t}=0,$ $ {\bf t}^\top A_1{\bf t}=0$ for all ${\bf t} \in T_{im},$
 and ${\bf e}^\top_1A_0{\bftau}(1)=0,$ $ {\bf e}^\top_1 A_1{\bftau}(1)=1.5.$  Hence,
${\bf a}({\bf t})={\bf 0}\ \forall {\bf t}  \in T_{im}$ and ${\bf b}(k_0,j_0)=(0,\,  {1.5})^\top$ for $k_0=1$, $j_0=1$, $k_0\in M(j_0).$  Thus we obtain
$${\bf b} (k_0,j_0)\not \in {\rm cone } \{ {\bf a} ({\bf t} ), {\bf t}  \in T_{im}\},$$
wherefrom we conclude that condition {\bf I}) does not hold true and, consequently,  {the} system
under consideration  does  not yield  the   uniform LP duality.

Let show this directly.  Since  the  system  $\ca({\bf x} )\in \cop$  with data (\ref{80}) has a unique feasible solution  ${\bf x } = x_{1}=0$,
then the corresponding primal problem ({\bf P}) has the optimal solution ${ \bf x }^* = x_{1}^*=0$ with $Val({\rm {\bf P}})=0$ for any objective function $ {\textbf c^{{\top}}
\textbf x^* =}c_{1}x_{1}^*{, \ c_1\in \mathbb R}$.  The corresponding dual problem ({\bf P}) takes the form
\be \max(-A_0\bullet U)\ \mbox{ s.t. } A_1\bullet U=c_1,\; U\in \cp.\label{D-1}\ee
 Suppose that  the  system $\ca({\bf x})\in \cop$ yields  the  uniform LP duality. Hence  the  dual problem should have an optimal solution $U^0$  such that
 \be U^0=\sum\limits_{i \in I}\alpha_i{\bf t}(i)({\bf t}(i))^\top\ \mbox{ with } \alpha_i>0,\, {\bf t}(i)\in T, \, i \in I,\ -A_0\bullet U^0=0,\ A_1\bullet
 U^0=c_1.\label{81}\ee
 Since ${\bf t}^\top A_0{\bf t}=a(t_1-t_3)^2$ for all $ {\bf t} \in \mathbb R^3,$ we conclude that  the
equality $ -A_0\bullet U^0=0$ implies  the equalities
 $ t_1(i)=t_3(i) \mbox{ for all } i \in I$  and hence  $({\bf t}(i))^\top A_1{\bf t}(i)=0$ for all $ i\in I$. Then  $ A_1\bullet U^0=0.$
 Thus we have  shown that for any $c_{1}\not =0$,  the dual problem has no solutions satisfying
relations (\ref{81})  which permits to  conclude that the  system  $\ca({\bf x})\in \cop$  with  the  data  defined in  (\ref{80}) does not yield  the   uniform LP duality.

\vspace{3mm}

 Note that in \cite{Pataki2020-1} (see page 3), it is stated that  for  {the}  SDP
systems $\ca({\bf x})\in {\cal S}^p_+$  with $n=1$, the  uniform   LP  duality  property is always satisfied.
  In our example,   for  $n=1$  we present  the CoP system $\ca({\bf x})\in \cop $  that does not yield the   uniform
	  LP  duality. This   further confirms the fact that      CoP systems are much  more  complex  (more {\it pathological}) than
	  SDP  systems.

\vspace{3mm}

{\bf Example 3.}  Consider a system $\ca({\bf x})\in \cop$ with the following data:
\be n=1, \ p=3, \ A_0=
\left(\begin{array}{ccc}
a& 0 &-a\cr
0 & 0& 0\cr
-a& 0 &a\end{array}\right),  \
A_1=
\left(\begin{array}{ccc}
1& 1 & -1\cr
1 & 0& 1\cr
-1& 1 &1 \end{array}\right)\ \mbox {with }a>0.\label{84}\ee
This system admits  feasible solutions   ${\bf{x}} = x_{1}\geq 0.$ It is easy to check that the set  $T_{im}=\{{\bftau}(j), j \in J\}$ consists of two vectors
 ${\bftau}(1)=0.5(1,\, 0,\, 1)^\top$, ${\bftau}(2)=(0,\, 1,\, 0)^\top$, and   that here  $J=\{1,2\}.$
 Since
 $\ca({\bf x}){\bftau}(1)=(0,\, x_{1},\, 0)^\top$ and $\ca({\bf x}){\bftau}(2)=(x_{1},\, 0,\, x_{1})^\top$, we conclude that
 the sets $M(j)$, $N(j),$ $j \in J=\{1,2\},$  take the form $M(1)=\{1,\,3\}$, $N(1)=\{2\}$, $M(2)=\{2\}$, $N(2)= \{1,\,3\} $,
  and   the   vectors ${\bf x}^*(k,j)\in X$ and sets $T_{im}(k,j)$, $ k \in N(j),$ $ j\in J,$  satisfying (\ref{W18}) are as follows:
$${\bf x}^*(k,j)=0, \ T_{im}(k,j)=\{{\bf t} \in T: t_1=t_3\}   \ \forall k \in N(j),\, j \in J.$$
  It is evident that the system with data (\ref{84}) does not satisfy  the  Slater condition.

 It follows from the equalities $A_0{\bftau}(1)=(0,\, 0,\, 0)^\top$, $A_1{\bftau}(1)=(0,\, 1,\, 0)^\top$,
$A_0{\bftau}(2)=(0,\, 0,\, 0)^\top$ and $A_1{\bftau}(2)=(1,\, 0,\, 1)^\top$,
that
$$ {\bf b}(k,j)=\left(\begin{array}{c}
0\cr
0\end{array}\right) \ \forall k \in M(j), \, j \in J;\
 {\bf b}(k,j)=\left(\begin{array}{c}
0\cr
1\end{array}\right), \ \forall k \in N(j), \, j \in J.
$$

Let us check  the conditions I) and II).   Here  we have
$$ {\bf b}(k,j)=\left(\!\begin{array}{c}
0\cr
0\end{array}\!\right)\!\in {\rm cone}\{{\bf a}({\bf t}),\,{\bf t}\in  T_{im}\}\!=\! {\rm cone}\{{\bf a}({\bftau}(j))=\left(\!\begin{array}{c}
0\cr
0\end{array}\!\right),\, j \in J\}\, \forall k \in M(j), j \in J,$$
\bea&{\bf b}(k,j)\!=\!\left(\begin{array}{c}
0\cr
1\end{array}\right)\in {\rm cone}\Bigl\{\left(\begin{array}{c}
0\cr
4t_2t_3\end{array}\right), {\bf t}\in   T_{im}(k,j) \Bigl\}
 \ \forall k \in N(j), \, j \in J.\nonumber\eea
 Thus,  we  can conclude that the  conditions I) and II)   hold true and, consequently, the  system under consideration yields  the
  uniform  {LP} duality despite it does not satisfy  the  Slater condition.

Let us check  the uniform  LP  duality  directly.  Since the system  $\ca({\bf x})\in \cop$  with the  data  defined in  (\ref{84}), has the set of feasible solutions in the
form $X=\{ x_1\geq 0\}$, then the corresponding primal problem ({\bf P}) has the finite  optimal value   $Val({\bf P})=0$ only for the objective function $c_1x_1$ with
$c_1\geq 0.$  The corresponding dual problem ({\bf D}) takes the form (\ref{D-1}).

 For a given $c_1\geq 0$,  set $U^0(c_1)=\frac{c_1}{4}\bar{\bf t}\,\bar{\bf t}^\top$ with $\bar {\bf t}= (1,\,1,\,1)^\top.$
It is easy to check that
 $$ U^0(c_1)\in \cp,\ -A_0\bullet U^0(c_1)=0,\ A_1\bullet U^0(c_1)=c_1.$$
Hence, for any $c_1\geq 0$,  $U^0(c_1)$  is a feasible solution of the corresponding dual problem ({\bf D}) and $Val({\bf P})=Val({\bf D}).$
 Thus, in fact, system  $\ca({\bf x})\in \cop$  with  data (\ref{84}) yields the uniform LP duality.

\section{ On a relationship  of the obtained results  with  the uniform duality for  SIP}\label{s:6}
 Consider a general linear SIP problem in the form
$$ {\bf P}^*_{SIP}: \qquad \qquad \qquad \displaystyle\min_{x \in \mathbb R^n } \ {\bf c}^{\top}{\bf x} \;\;\;
\mbox{s.t. } (1,{\bf x}^\top){\bf a}({\bf t})\geq 0 \ \forall {\bf t} \in T, \qquad \qquad $$
 where $T$ is an index set  and  ${\bf a}({\bf t})=(a_m({ \bf t}), m=0,1,...,n)^\top,$ $ {\bf t} \in T.$
Denote
$$G:=\left\{ {\bf a}({\bf t}), {\bf t} \in T; (1,{\bf 0}^\top_n)^\top\right\}.
$$

The following theorem is proved in   \cite{Duffin}  (see  Theorem 3.2,  conditions (ii) and (iv)).

\begin{theorem}\label{T-Duf}
The  consistent  constraint system of   the  problem  (${\bf P}^*_{SIP}$)
 yields  the  uniform LP duality iff
 \be {\rm cone }(G)={\rm cone} (F\cup W)\label{09-12-4}\ee
 with some $F\subset \mathbb R^{n+1}$ and $W\subset \mathbb R^{n+1}$  satisfying the following conditions:  \begin{itemize} \item $F$ is finite and ${\rm cone}(F)$ is
a linear space which is also contained in ${\rm cone}\, G$,
\item  $W$ is compact, and there exists a vector $\bar{\bf  x}\in \mathbb R^n$ such that
\be \begin{split} & s_{0}+{\bf s}^\top \bar{\bf x}=0 \ \forall  (s_{0}, {\bf s}^\top)^\top \in F, \ s_{0} \in \mathbb R,\\
& t_{0}+{\bf t}^\top \bar{\bf x} >0 \ \forall   (t_{0}, {\bf t}^\top)^\top \in W,\  t_{0}\in \mathbb  R.\end{split}\label{09-12-5}\ee
\end{itemize}
Moreover, whenever $(s_{0}, {\bf s}^\top)^\top \in F$, the equality ${\bf s}^\top {\bf x} =- s_{0} $ is implied by the
constraint system of  the problem   (${\bf P}^*_{SIP}$)  and the equality ${\bf s}^\top  {\bf x} = 0 $ is implied by the
homogeneous system
$(0,{\bf x}^\top){\bf a}({\bf t})\geq 0, {\bf t} \in T$. Furthermore,  it holds:
\be {\bf s}^\top {\bf t} + s_{0}t_{0} = 0 \mbox{ if } (s_{0},{\bf s}^\top)^\top \in F \mbox{ and } ( t_{0}, {\bf t}^\top)^\top \in W.\label{09-12-1}\ee
\end{theorem}

 It follows from the equivalent description (\ref{copT}) of the cone $\cop$ that  the problem  ({\bf P}) is equivalent to the   linear    SIP problem  (${\bf P}^*_{SIP}$)
where  the set  $T$ and  the   vector ${\bf a}({\bf t})$     are  defined in (\ref{SetT}) and (\ref{a-t}), respectively. Let us denote this special SIP problem by (${\bf
P}_{SIP}$).

Since  the  problem (${\bf P}_{SIP}$)   is a   special  case of a general linear SIP problem,
 the statements of  Theorem \ref{T-Duf} should be   satisfied  for  the  problem (${\bf P}_{SIP}$)  as well.

In Theorem \ref{T-Duf},  nothing  is  said about how to build the sets $F$ and $W$ mentioned in the theorem.
 Obviously, it is interesting to know  how to  find  these sets for our CoP problem. The following theorem gives an answer to this question  for the  CoP problem  under
 consideration.

 \begin{theorem}
 Given the   problem (${\bf P}_{SIP}$),
  the  sets $F$ and $W$, mentioned in Theorem \ref{T-Duf}, can  be chosen as follows:
\bea& F:=\{ {\bf b}(k,j), k \in M(j),j\in J\}, \ W:={\rm Pr}(\widetilde{W},L^\perp), \label{09-12-2}\\
&\widetilde{W}:=\{{\bf a}( {\bf t}),{\bf  t} \in \Omega\}\cup\{{\bf b}(k,j), k \in N(j),j \in J\}\cup (1,{\bf 0}^\top_n)^\top,\label{09-12-3}\eea
where $L:={\rm span} \,F$,   $L^\perp$ is the orthogonal complement to $L$,  the set $\Omega$ is defined in (\ref{Omega}),  and   ${\rm Pr}(\widetilde{W},L^\perp)$ is the
projection of the set $\widetilde{W}$ onto  $L^\perp$.
\end{theorem}
{\bf Proof}. First, let us show  that the consistent  constraints  system of   the  problem   (${\bf P}_{SIP}$)  yields  the  uniform LP duality iff
\be  {\rm cone} (F\cup\widetilde{ W})= {\rm cone }(G)\label{zzz}\ee
 with the sets $F$ and $\widetilde{W}$ defined in (\ref{09-12-2}) and (\ref{09-12-3}).

 It was shown above that the problem (${\bf P}_{SIP}$) is equivalent to the problem (${\bf P}_*$) (see (\ref{W6})) that can be rewritten in the form
 $$ {\bf P}_*: \ \ \displaystyle\min_{{\bf x} \in \mathbb R^n } \ {\bf c}^{\top}{\bf x} \;\;\;
\mbox{s.t. } (1,{\bf x}^\top){\bf a}({\bf t})\geq 0 \ \forall {\bf t} \in \Omega,\ (1, {\bf x}^\top){\bf b}(k,j)\geq 0, k \in N_*(j), j \in J.$$
Notice that the problem (${\bf P}_*$) yields  the   uniform  LP  duality, the problems (${\bf P}_{SIP}$) and (${\bf P}_*$)
 have the same  feasible  sets
 and, consequently, the same optimal values of the cost functions.

 For a  fixed ${\bf t}\in T$, consider  the consistent  problem (${\bf P}_*$) with  ${\bf c}^\top=(c_m={\bf t}^\top A_m{\bf t},\, m=1,...,n).$
Then  for all ${\bf x}\in X$,
$${\bf c}^\top {\bf x}=\sum\limits_{m=1}^n {\bf t}^\top A_m {\bf t}x_{m}\geq - {\bf t}^\top A_0 {\bf t}>-\infty. $$
Consequently
  $Val({\bf P}_{SIP})=\beta - {\bf t}^\top A_0 {\bf t}$ with some $\beta\geq 0$. Taking into account that  the
	problem (${\bf P}_*$) yields  the   uniform LP duality, we conclude that
  there exist numbers and vector
$$\alpha_i,\ {\bf t}(i)\in \Omega,\ i \in I,   \ \lambda_k(j), \, k\in N_{*}(j),\ j \in J,$$
such that $\alpha_i> 0,i \in I,$ $|I|<\infty,$ $ \lambda_k(j)\geq 0, \, k\in N_{*}(j),\ j \in J,$ and equality (\ref{W12}) holds true
with
   ${\bf c}=(c_m={\bf e}^\top_m{\bf a}({\bf t}),\, m=1,...,n)$ and $Val({\bf P}_{SIP})=\beta- {\bf t}^\top A_0 {\bf t}$.
   This implies that
   $${\bf a}({\bf t})=\sum\limits_{i\in I}\alpha_i{\bf a}({\bf t}(i))+ \sum\limits_{j\in J} \sum\limits_{k\in
    N_{*}(j)}\lambda_k(j){\bf b}(k,j)+(1,{\bf 0}^\top_n)^\top \beta\in {\rm cone} (F\cup\widetilde{ W}).$$
  Since  the inclusion  above is satisfied for all ${\bf t}\in T$, we conclude
 that for  the   consistent  problem (${\bf P}_{SIP}$),  it holds:
\be {\rm cone }(G)\subset{\rm cone} (F\cup\widetilde{ W}).\label{09-12-17}\ee

 \vspace{3mm}

On the other hand, due to Proposition \ref{P1}, the consistent system of   the  problem   (${\bf P}_{SIP}$)  yields  the    uniform  LP  duality iff
$${\bf b}(k,j)\in {\rm cone}\{{\bf a}({\bf t}), {\bf t}\in T\}\subset {\rm cone} (G)\  \forall \,k \in N_*(j),\  \forall\, j\in J.$$
Taking into account these inclusions and  the  definitions of the sets $G$, $F$ and $\widetilde{W},$ we obtain  that
the consistent system of   the   problem   (${\bf P}_{SIP}$)  yields  the   uniform LP duality iff
$$ {\rm cone} (F\cup\widetilde{ W})\subset {\rm cone }(G).$$
Now taking into account inclusion (\ref{09-12-17}), we conclude that the consistent system of  the   problem   (${\bf P}_{SIP}$)  yields  the  uniform LP duality iff
condition (\ref{zzz})  is satisfied
 with the sets $F$ and $\widetilde{W}$ defined in (\ref{09-12-2}) and (\ref{09-12-3}).

\vspace{3mm}

 By construction, the set $F$ is finite and the set $\widetilde{W}$ is compact.
 Relations (\ref{09-12-5}) hold true  with $\bar {\bf x}={\bf x}^*$ where ${\bf x}^*$ is defined in (\ref{W7}).
It follows from Proposition \ref{P0} that ${\rm cone} \,F={\rm span} \,F$ and then it is evident that  the set  $ {\cal L}:={\rm cone} \,F$ is a subspace.

 Hence we have shown
 that the sets $F$ and $\widetilde{W}$ defined in (\ref{09-12-2}) and (\ref{09-12-3}) satisfy all statements of Theorem
\ref{T-Duf} except for the condition (\ref{09-12-1}).
 Taking into account that  the set ${\cal L}$ introduced above  is a subspace, it is easy to see that
to satisfy the condition (\ref{09-12-1}) it is enough to replace the set $\widetilde{W}$ by its projection onto ${\cal L}^\perp.$
 The theorem is proved.$\quad \Box$

\section{ Conclusions }\label{s:7}

 The main result of the paper is to establish the  necessary and sufficient conditions
that guarantee  the    uniform LP duality for linear CoP problems.
These conditions are obtained using the concept of immobile indices and the sets generated by them and are formulated in various equivalent
forms thereby expanding the scope of their application.  The examples illustrate how the conditions obtained can be applied
to confirm or deny the  {uniform} LP duality of a CoP system. The relationship between the uniform LP duality properties  for the related problems of CoP and SIP  is shown.

\appendix  
\section{Appendix}

\begin{proposition}\label{P01} The equalities (\ref{cond-bb-1})  and (\ref{W30})  are equivalent.
\end{proposition}
{\bf Proof.} Consider any ${\bf t}\in T_{im}$. From (\ref{W32}), it follows  that ${\bf t}\in T_{im}(s)$ for some $s\in S$  and ${\bf t}$ admits a representation (\ref{W20}).
Consequently,
$${\bf t}^\top\cb({\bf z} ){{\bf t}}=\sum\limits_{i\in J(s)}\sum\limits_{j\in J(s)}\alpha_i\alpha_j({\bftau}(i))^\top\cb({\bf z}){\bftau}(j).$$
Taking into account this equality, we conclude that  equalities (\ref{cond-bb-1}) imply equalities (\ref{W30}).

Now we suppose that   equalities  (\ref{W30}) hold true.
Since for any $s\in S$ and  for all $i\in J(s)$ and all $j\in J(s)$,  we have  $0.5({\bftau}(i)+{\bftau}(j))\in T_{im}$, it follows from equalities  (\ref{W30}) that $$
({\bftau}(i)+{\bftau}(j))^\top \cb({\bf z})({\bftau}(i)+{\bftau}(j))=0 \  \forall i\in J(s), \ \forall j\in J(s).$$
It is evident that these equalities imply (\ref{cond-bb-1}).
 $\ \Box$

\begin{proposition}\label{P02} Consider a  set $L:={\rm cone}\{{\bf a}(i),i\in {\cal I}\}$ where ${\bf a}(i)\in \mathbb R^s,$ $  {\cal I}$ is a set of indices (it is possible
that $| {\cal I}|=\infty$). The set $L$   is a subspace iff there exist a finite subset $\{{\bf a}(j), j \in J\} \subset \{{\bf a}(i), i \in  {\cal I}\},$ $|J|<\infty,$ and
numbers $\alpha_j,$ $ j\in J,$ such that
\be \alpha_j>0\ \forall j \in J, \ \sum\limits_{j\in J}\alpha_j{\bf a}(j)=0, \ {\rm rank}({\bf a}(j), j \in J)= r_*:={\rm rank}({\bf a}(i), i \in  {\cal I}).\label{1}\ee
\end{proposition}
{\bf Proof.} Suppose that there exist numbers   and vectors  $\alpha_j,$  ${\bf a}(j),$  $ j\in J,$ such that (\ref{1}) holds true. To show that $L$ is a subspace, we have to
show that $-{\bf d}\in L$ for any ${\bf d}\in L$.

Since ${\bf d} \in L$ and $ {\rm rank}({\bf a}(j), j \in J)=r_*$,  then ${\bf d}=\sum\limits_{j\in J}\beta_j{\bf a}(j)$ with some $\beta_j, $ $j\in J.$  Hence  $-{\bf
d}=-\sum\limits_{j\in J}\beta_j{\bf a}(j)$ and taking into account (\ref{1}) we get
$$-{\bf d}=-\sum\limits_{j\in J}\beta_j{\bf a}(j)=\sum\limits_{j\in J}(M\alpha_j-\beta_j){\bf a}(j)=\sum\limits_{j\in J}\bar\alpha_j{\bf a}(j),$$
where $M:=\max\{\beta_j/\alpha_j , j\in J, 0\}$,
$\bar\alpha_j:=M\alpha_j-\beta_j\geq 0$ $\forall\, j \in J.$  Then  $-{\bf d}\in L.$
Thus we have shown that if relations (\ref{1}) hold true then $L$ is a subspace.

Now,  assuming  that  $L$ is a subspace, let us show that relations (\ref{1})  hold true.

 It is evident that there exists a subset $ {\cal I}_*\subset  {\cal I} $ such that
$${\rm rank}({\bf a}(i), i \in  {\cal I}_*)=| {\cal I}_*|=r_*.$$

 Since $L$ is a subspace then $-{\bf a}(i)\in L$ for all $ i\in  {\cal I}_*.$  Hence, for   any  $i \in  {\cal I}_*$, there exist   a set $ {\cal I}(i)\subset  {\cal I}$ and
 numbers $\alpha_{ij}, j \in  {\cal I}(i),$ such that
$$-{\bf a}(i)=\sum\limits_{j \in  {\cal I}(i)}\alpha_{ij}{\bf a}(j),\ \alpha_{ij}>0, j \in  {\cal I}(i),\ | {\cal I}(i)|\leq r_*.$$
 Consider the following set of vectors:
 \be \{{\bf a}(i),\ {\bf a}(j),\, j \in  {\cal I}(i);\ i \in  {\cal I}_*\}.\label{f-set}\ee
 This set consists in a finite number of elements and, by construction,
 $${\rm rank}({\bf a}(i),\ {\bf a}(j),\, j \in  {\cal I}(i);\ i \in  {\cal I}_*)=r_*,$$
 $$\sum\limits_{j \in  {\cal I}_*}({\bf a}(i)+\sum\limits_{j \in  {\cal I}(i)}\alpha_{ij}{\bf a}(j))=\sum\limits_{j \in  {\cal I}_*}({\bf a}(i)-{\bf a}(i))=0.$$
 \ Thus  we get that  relations (\ref{1})  hold true with the finite set of vectors (\ref{f-set}).
  $\ \Box$

\begin{proposition}\label{P03} Consider a  set $\{{\bf a}(i),i\in  {\cal I}\}$ where ${\bf a}(i)\in \mathbb R^s,$ $  {\cal I}$ is a set of indices (it is possible that $|
{\cal I}|=\infty$). Suppose that  there exist a finite subset $\{{\bf a}(j), j \in J\} \subset \{{\bf a}(i), i \in  {\cal I}\},$ $|J|<\infty,$ and numbers $\alpha_j,$ $ j\in
J,$ such that relations (\ref{1}) are satisfied. Then for any $\Delta J\subset   {\cal I},$ $|\Delta J|<\infty,$ there exist numbers $ \bar \alpha_j,\,  j \in J, $ such that
\be \bar \alpha_j>0\ \forall j \in J, \ \sum\limits_{j\in J}\bar \alpha_j{\bf a}(j)+\sum\limits_{q\in \Delta J}{\bf a}(q)=0.\label{1W}\ee
\end{proposition}
{\bf Proof.} It follows from (\ref{1}) that for all $q\in \Delta J,$ there exist numbers $\beta_j(q),\, j\in J,$ such that
$$ {\bf a}(q)=\sum\limits_{j\in J}\beta_j(q){\bf a}(j).$$
 Then, taking into account the equality in (\ref{1}), we conclude  that for any $M\in \mathbb R,$ the  equality
\be \sum\limits_{j\in J}( M\alpha_j-\sum\limits_{q\in \Delta J}\beta_j(q) ){\bf a}(j)+\sum\limits_{q\in \Delta J}{\bf a}(q)=0\label{2W}\ee
 holds true.
Let us denote $\Delta \alpha(j)=\sum\limits_{q\in \Delta J}\beta_j(q)$, $ j \in J,$  and  choose $M:=\max\{\Delta \alpha(j)/ \alpha_{j},j \in J, 0\}+1.$
Then (\ref{2W}) implies (\ref{1W}) with $\bar \alpha_j=M \alpha_j-\Delta \alpha(j)\geq  \alpha_j>0, j \in J.$  $\ \Box$
\begin{proposition}\label{P04} The following  equality  holds true:
$$ r_{im}:={\rm rank } ({\bf a}({\bf t}), {\bf t} \in T_{im})={\rm rank } ({\bf a}({\bftau}(l)+{\bftau}(q)), (l,q)\in V_0)=:\bar r,$$ 
where $V_0$ is defined in (\ref{V0}).
\end{proposition}
{\bf Proof. } Since $\frac{1}{2}({\bftau}(l)+{\bftau}(q))\in T_{im}$  for all   $(l,q)\in V_0,$ it is evident that $\bar r\leq r_{im}.$

Suppose that $\bar r< r_{im}.$ Then there exists $ {\bf z}\in \mathbb R^{n+1}$ and ${\bf t}\in T_{im}$  such that $ {\bf z}^\top {\bf a}({\bf t})\not=0$ and
$${\bf z}^\top  {\bf a}({\bftau}(l)+{\bftau}(q))=0\ \forall (l,q)\in V_0.$$
Notice that the latter equalities imply the equalities
\be \sum\limits_{m=0}^nz_m({\bftau}(l))^\top A_m{\bftau}(q)=0\ \ \forall (l,q)\in V_{0}. \label{W43}\ee

Since ${\bf t}\in T_{im}$,  then ${\bf t}\in T_{im}(s)$ with some $s\in S$ and ${\bf t}$ admits a representation (\ref{W20}).
Therefore, with respect to (\ref{W43}),  we conclude that
\bea &0\not={\bf z}^\top {\bf a}({\bf t})=\sum\limits_{m=0}^nz_m{\bf t}^\top A_m{\bf t}=
\sum\limits_{m=0}^nz_m\Bigl(\sum\limits_{l\in J(s)}\alpha_l{\bftau}(l)\Bigl)^\top A_m\Bigl(\sum\limits_{q\in J(s)}\alpha_q{\bftau}(q)\Bigl)\nonumber\\
&=\sum\limits_{l\in J(s)}\sum\limits_{q\in J(s)}\alpha_l \alpha_q \sum\limits_{m=0}^nz_m({\bftau}(l))^\top A_m{\bftau}(q)=0.\nonumber\eea
  The  contradiction  obtained  shows that   $\bar r= r_{im}.$  $\ \Box$

\begin{proposition}\label{P05} { For any ${\bf t}\in T_{im}$ there exist numbers $\beta_{lq}  \in \mathbb R,$  where  $(l,q)\in V_0,$  such that
\be {\bf t}{\bf t}^\top=\sum\limits_{(l,q)\in V_0}\beta_{lq}({\bftau}(l)+{\bftau}(q))({\bftau}(l)+{\bftau}(q))^\top.\label{W51}\ee}
\end{proposition}
{\bf Proof.} Consider any ${\bf t}\in T_{im}$. Then ${\bf t}\in T_{im}(s)$ with some $s\in S$ and ${\bf t}$ admits a representation (\ref{W20}). Consequently
\be\begin{split}&{\bf t}{\bf t}^\top=\bigl(\sum\limits_{l\in J(s)}\alpha_l\bftau(l)\bigl)\bigl(\sum\limits_{q\in J(s)}\alpha_q\bftau(q)\bigl)^\top=
\sum\limits_{l\in J(s)}\sum\limits_{q\in J(s)}\alpha_l\alpha_q\bftau(l)\bftau^\top(q)\\
&\qquad =\sum\limits_{l\in J(s)}\alpha^2_l\bftau(l)\bftau^\top(l)+\sum\limits_{(l,q)\in V(s), l\not=q}\alpha_l\alpha_q
 [\bftau(l)\bftau^\top(q)+\bftau(q)\bftau^\top(l) ].\end{split}\label{12-12-111}\ee
It is evident  that for any $(l,q)\in V(s),$ $l\not=q,$ we have
$$\bftau(l)\bftau^\top(q)+\bftau(q)\bftau^\top(l)=
({\bftau}(l)+{\bftau}(q))({\bftau}(l)+{\bftau}(q))^\top-{\bftau}(l){\bftau}^\top(l)-{\bftau}(q){\bftau}^\top(q).$$
These equalities together with (\ref{12-12-111}) imply the equality (\ref{W51}). $\ \Box$

\end{document}